\documentclass[11pt]{article}

\usepackage[english]{babel}

\usepackage[letterpaper,top=2cm,bottom=2cm,left=3cm,right=3cm,marginparwidth=1.75cm]{geometry}

\usepackage{amsmath}
\usepackage{amssymb}
\usepackage{graphicx}
\usepackage[colorlinks=true, allcolors=blue]{hyperref}
\usepackage[dvipsnames]{xcolor}
\usepackage{todonotes,setspace}
\usepackage{tikz}
\usepackage{bbm}
\usepackage{cite}
\usepackage{subcaption}
\usepackage{centernot}

\newcounter{mycomment}
\newcommand{\mycomment}[2][]{%
	\refstepcounter{mycomment}%
	\ifcase\numexpr \themycomment- 4*((\themycomment+3)/4 -1) \relax
	{%
		\setstretch{0.7}
		\todo[color=SeaGreen,size=\scriptsize]{%
			\textbf{[\uppercase{#1}\themycomment]:} #2}%
	}
	\or
	{%
		\setstretch{0.7}
		\todo[color=Salmon,size=\scriptsize]{%
			\textbf{[\uppercase{#1}\themycomment]:} #2}%
	}
	\or
	{%
		\setstretch{0.7}
		\todo[color=Thistle,size=\scriptsize]{%
			\textbf{[\uppercase{#1}\themycomment]:} #2}%
	}
	\else
	{%
		\setstretch{0.7}
		\todo[color=Tan,size=\scriptsize]{%
			\textbf{[\uppercase{#1}\themycomment]:} #2}%
	}
	\fi
}

\usepackage{amsthm}

\theoremstyle{definition}
\newtheorem{exmp}{Example}[section]

\newtheorem{corollary}{Corollary}[section]
\newtheorem{theorem}{Theorem}[section]
\newtheorem{proposition}{Proposition}[section]

\newtheorem{lemma}{Lemma}[section]
\usepackage{mathtools}
\usepackage{tkz-euclide}

\usepackage{orcidlink}

\usepackage{amsmath}
\usepackage{amsfonts}
\usepackage{bbm}
\usepackage{pgfplots}
\pgfplotsset{compat=1.18}
\usepackage{tikz}
\usepackage{enumerate}
\usepackage{graphicx} 

\title{Analytic properties of distributions associated with multi-state  Random Number Generators} 
\author{Douglas T. Pfeffer\,\orcidlink{0000-0001-8287-2686} and Volodymyr Zabolotnyi}
\date{}

\begin{document}
	\maketitle
	
	\begin{abstract}
		We analyze the analytic properties of probability distributions associated with possibly unfair, multi-state Random Number Generators (i.e., dice). Specifically, given a $q$-sided die, if $x_i\in \{0,\ldots, q-1\}$  denotes the outcome of the $i^{\text{th}}$ toss, then $F(x)=\mathbb{P}[X\leq x]$ is the associated cumulative distribution function---a probabilistic generalization of the famed Cantor function. We first establish easily computable bounds for the $L^1$-distance between the possibly unfair $F$ and the uniform distribution. We then prove that $F$ is both upper- and lower-H\"{o}lder continuous with sharp H\"{o}lder exponents, which we apply to bound the Hausdorff dimension of the support of $dF$. Finally, we establish recursive formulae for the moments, $\mathbb{E}[X^n]$, of $F$ and use these moments to construct step-function envelopes that converge uniformly to $F$ from above and below.
	\end{abstract}
	
	\renewcommand{\thefootnote}{\fnsymbol{footnote}} 
	\footnotetext{\emph{2020 Mathematics Subject Classification.} Primary 60E05, Secondary 26A30,  26A16}     
	\renewcommand{\thefootnote}{\arabic{footnote}} 
	
	\renewcommand{\thefootnote}{\fnsymbol{footnote}} 
	\footnotetext{\emph{Key words.} Random Number Generators, Cantor function, singular distribution, H\"{o}lder continuity, moment generating function}     
	\renewcommand{\thefootnote}{\arabic{footnote}} 
	
	\newcommand{\Addresses}{{
			\bigskip
			\footnotesize
			
			\setlength{\parindent}{0cm}
			\noindent Douglas T. Pfeffer\,\orcidlink{0000-0001-8287-2686}\par\nopagebreak\textsc{Department of Mathematics, University of Tampa, 401 W. Kennedy Blvd. Tampa, FL 33606}\par\nopagebreak
			e-mail: \texttt{dpfeffer@ut.edu}
			
			\medskip
			
			Volodymyr Zabolotnyi \,\par\nopagebreak  \textsc{Department of Mathematics, University of Tampa, 401 W. Kennedy Blvd. Tampa, FL 33606
			}\par\nopagebreak
			e-mail: \texttt{volodymyr.zabolotnyi@spartans.ut.edu}
			
	}}
	
	\section{Introduction and Background}
	
	Pseudo-Random Number Generator (pRNG) are commonly used for quickly generating random numbers, where they typically use a \textit{seed} value to generate a sequence of numbers through iteration.  Whilst this sequence appears random, they are entirely determined by the seed and the pRNG.  Though useful for simple applications, the deterministic nature of these pRNGs can be a cause for concern. A number generator that produces truly random output is therefore desirable in fields such as computing, where probabilistic devices and true Random Number Generators (tRNGs) are an area of active study \cite{misra2022probabilistic,chowdhury2023accelerated,chowdhury2023full,aadit2022massively,kaiser2021probabilistic}. Certain devices, such as magnetic tunnel junctions \cite{rehm2022stochastic,liu2022random} or tunnel diodes \cite{bernardo2017extracting}, primarily operate as a Bernoulli random variable (i.e., a coin flip) with probability $p$ of returning 1 and a probability $1-p$ of returning 0. However, they can also operate as multi-state devices (i.e., $q$-sided dice), as seen in \cite[Fig. 2e]{leonard2023stochastic,leonard2022shape}. 
	Though the bulk of current investigations examine coin-like devices, it was argued in \cite{Pfeffer2025Dice} that dice-like devices might be more attractive for various reasons. In that work, the authors revisited classic results in the Bernoulli case for a class of singular measures and established $q$-sided dice analogues. Importantly, since human-made devices may exhibit small amounts of unfairness, the authors conducted their analysis while allowing the probability weights of each outcome to be possibly unfair (e.g., a theoretical coin is analyzed whilst allowing $p\neq 0.5$). We record some of their results presently.
	
	For $q\geq 2$, consider a series of Bernoulli outcomes $x_i$ on $\{0,1,\ldots, q-1\}$. In the independent and identically distributed (i.i.d.) case, let $p_j=\mathbb{P}[x_i=j]$ for $0\leq j\leq q-1$.  If $p_j=1/q$ for all $j$, we say we are rolling a `fair' die. Given $n$ outcomes, define
	\[
	X_n := \frac{x_1}{q^1} + \frac{x_2}{q^2} + \cdots + \frac{x_n}{q^n} = \sum_{i=1}^n \frac{x_i}{q^i}.
	\]
	Each $X_n$ is a decimal number in $[0,1]$ formed from a $q$-ary encoding of $n$ outcomes $\{x_i\}_{i=1}^n$.
	
\begin{figure*}[t]
		\centering
		
		\begin{subfigure}[b]{0.25\textwidth}
			\centering
			\includegraphics[width=\linewidth]{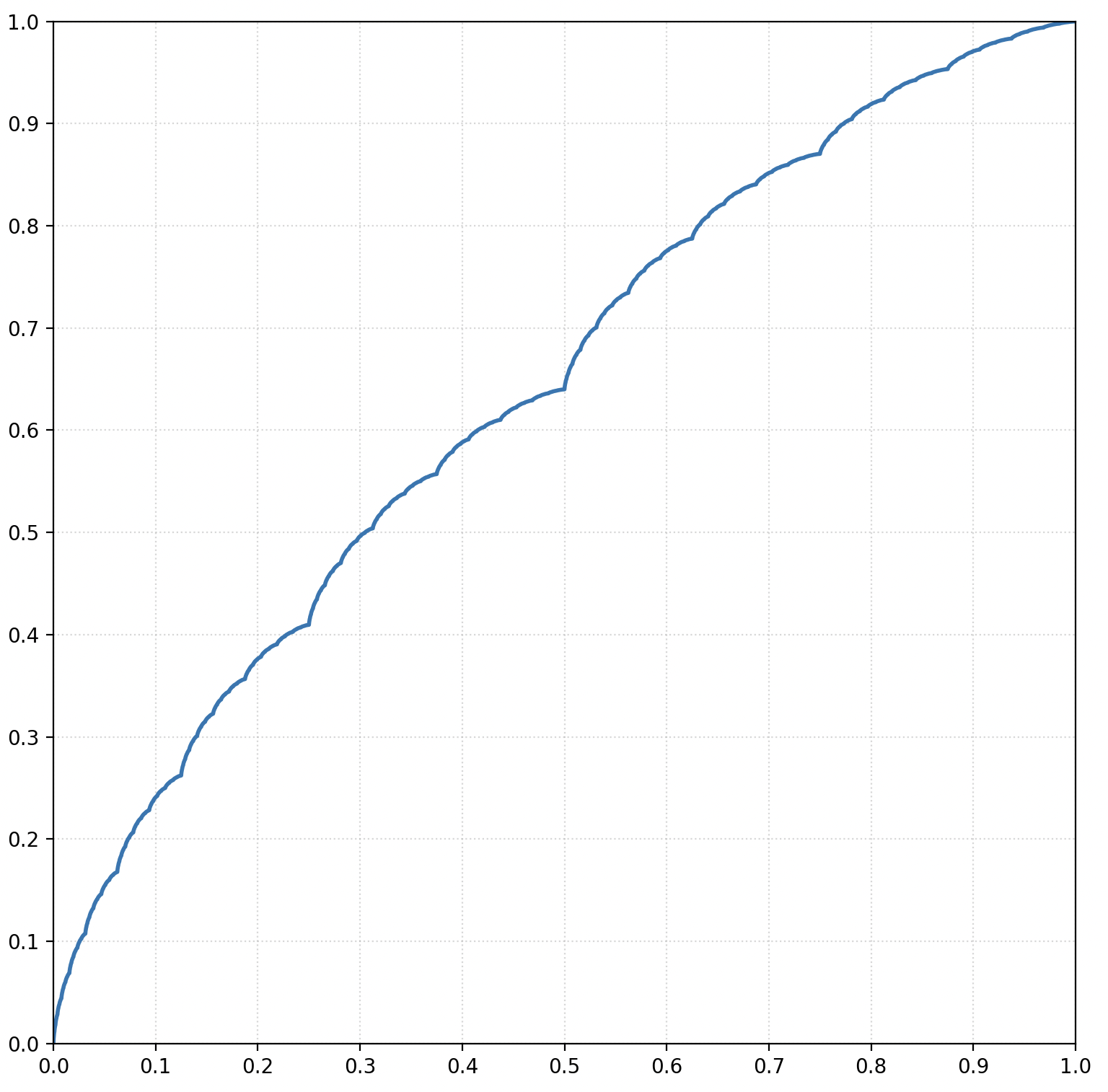}
			\caption{$P=\{0.64,\,0.36\}$}
			\label{fig:cdf_1}
		\end{subfigure}\hfill
		\begin{subfigure}[b]{0.25\textwidth}
			\centering
			\includegraphics[width=\linewidth]{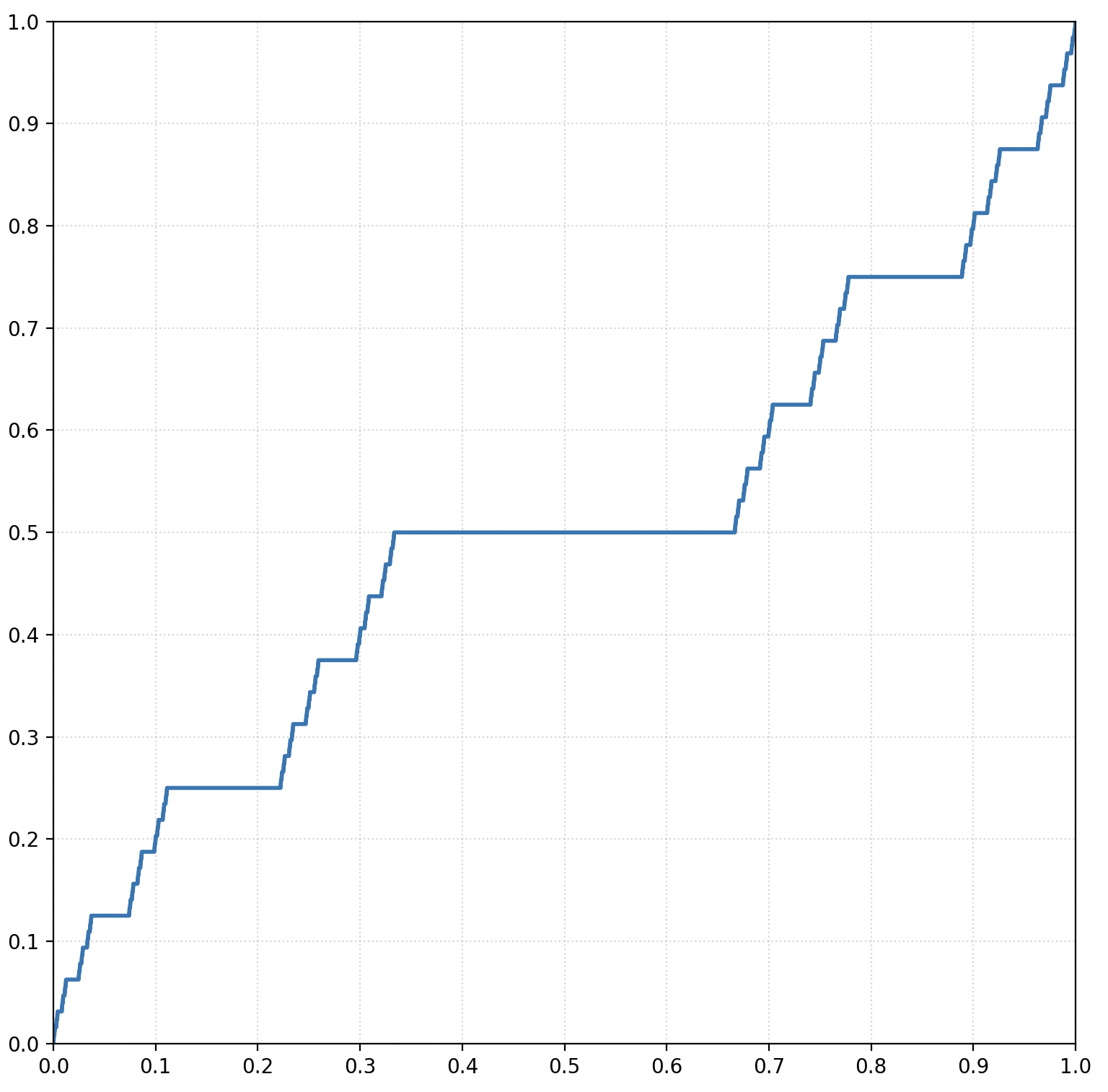}
			\caption{$P=\{0.50,\,0,\,0.50\}$}
			\label{fig:cdf_2}
		\end{subfigure}\hfill
		\begin{subfigure}[b]{0.25\textwidth}
			\centering
			\includegraphics[width=\linewidth]{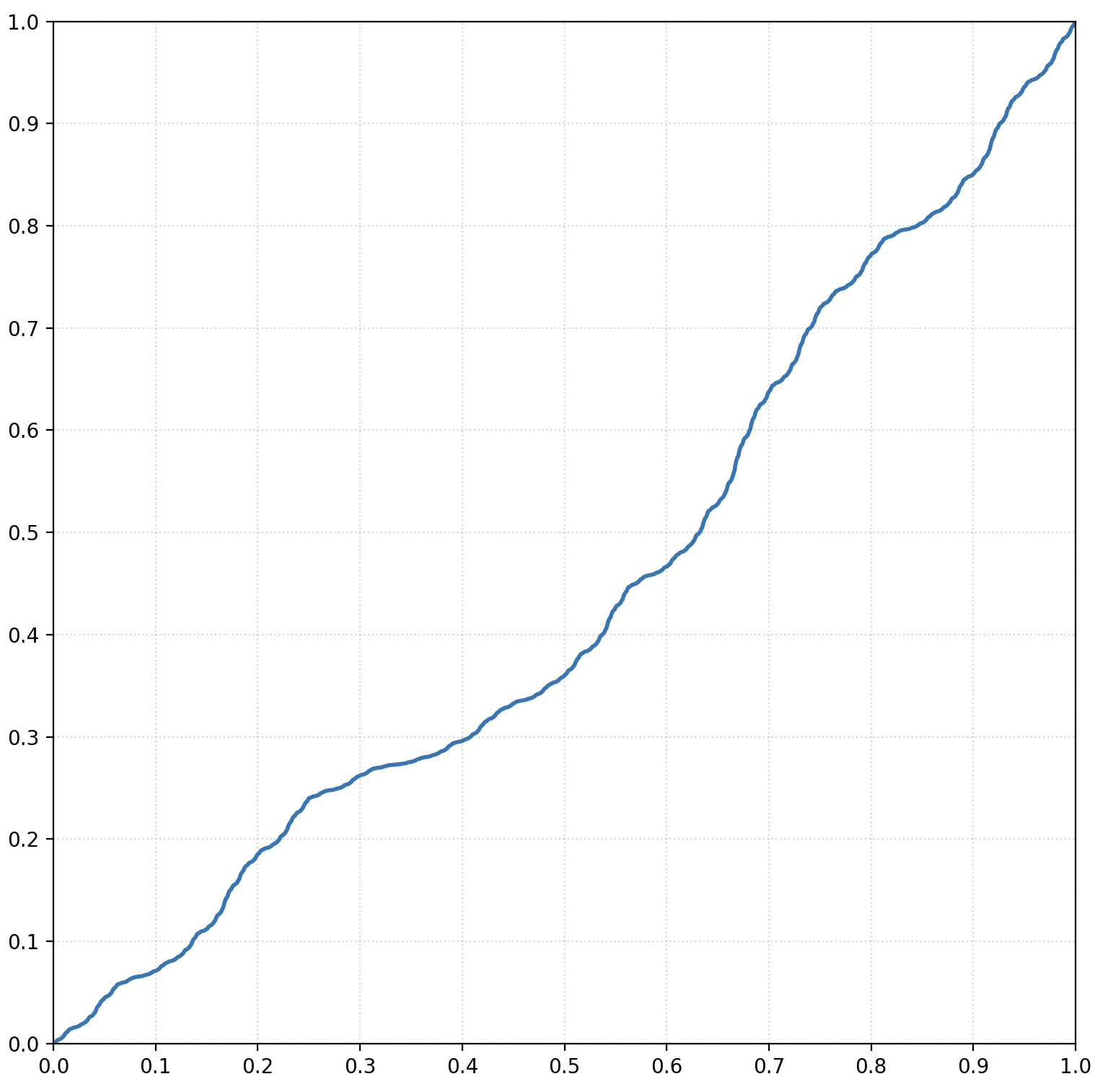}
			\caption{$P=\{0.24,\,0.12,\,0.36,\,0.28\}$}
			\label{fig:cdf_3}
		\end{subfigure}\hfill
		\begin{subfigure}[b]{0.25\textwidth}
			\centering
			\includegraphics[width=\linewidth]{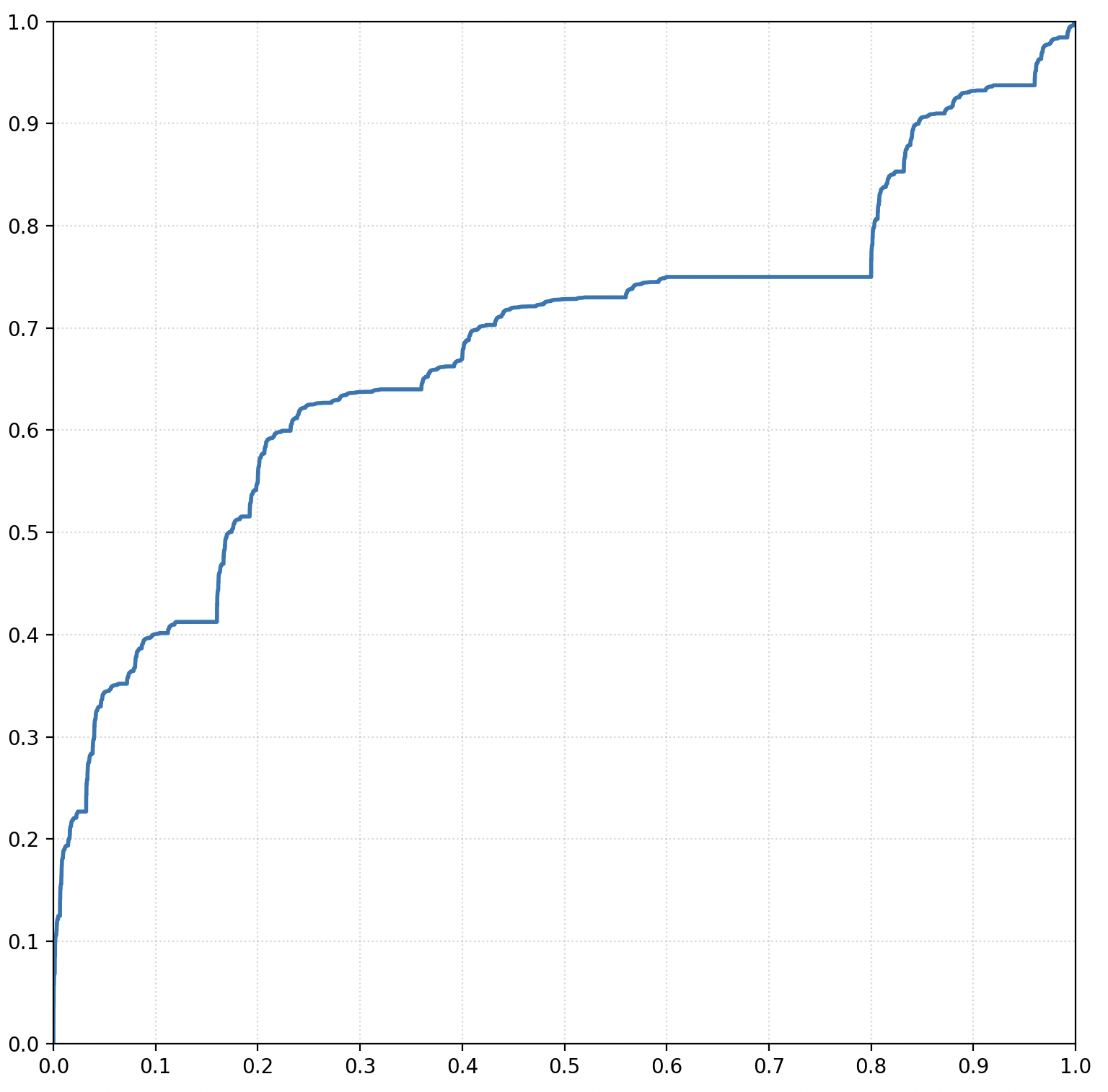}
			\caption{$P=\{0.55,\,0.12,\,0.08,\,0,\,0.25\}$}
			\label{fig:cdf_4}
		\end{subfigure}
		
		\caption{The cumulative distribution function $F(x)$ for four unfair dice, where the associated probability weights are given as $P=\{p_i\}$.}
		\label{fig:unfair_example}
	\end{figure*}

	Let $X := \lim_{n\to\infty} X_n$ be the encoded value in $[0,1]$ obtained by rolling our die infinitely many times. For $x\in [0,1]$, let  $F(x)=\mathbb{P}[X\leq x]$ denote the cumulative distribution function (CDF). As visual examples, we graph $F$ for four different die  
	in Figure \ref{fig:unfair_example}. The following result from \cite{Pfeffer2025Dice} helps characterize $F$:
	
	\begin{theorem}[Theorem 1 in \cite{Pfeffer2025Dice}]\label{maintheorem1}
		Consider a $q$-sided die, where $x_i\in\{0,\ldots,q-1\}$ denotes the outcome of the $i^\text{th}$ toss and $p_j:=\mathbb{P}[x_i=j]$ for $j\in\{0,\ldots,q-1\}$ with each $p_j<1$. Given $X = \sum_{i=1}^\infty x_iq^{-i}$, if $F(x)=\mathbb{P}[X\leq x]$ is the CDF obtained after tossing the die an infinite number of times, then
		\begin{enumerate}[(i)]
			\item $F$ is continuous;
			\item $F$ is strictly increasing on $[0,1]$ if and only if $p_j>0$ for all $j$; \label{item2thm1}
			\item $F'$ exists almost everywhere in $[0,1]$;
			\item If $p_j= \textstyle \frac{1}{q}$ for all $j$, then $F(x)=x$ on $[0,1]$ and
			\item If $p_k\neq \textstyle \frac{1}{q}$ for some $k$, then $F$ is singular.
		\end{enumerate}
		Additionally, $F$ is given by the following recursion formula:
		\begin{equation} \label{qsideddierecursionmain}
			F(x) = 
			\begin{cases} 
				p_0F(qx) &\text{if} \ 0\leq x \leq \textstyle \frac{1}{q}\\
				p_0+p_1F(qx-1) &\text{if} \ \textstyle \frac{1}{q} \leq x \leq \frac{2}{q}\\
				\ \ \ \ \ \ \ \ \vdots &\ \ \ \ \ \ \ \ \ \vdots\\
				(p_0+p_1+\ldots+p_{q-2}) + p_{q-1}F(qx-(q-1)) &\text{if} \ \textstyle \frac{q-1}{q} \leq x \leq 1\\
			\end{cases}.
		\end{equation}
	\end{theorem}

	\begin{exmp} \label{ex:cantor}
		If $p_0=p_2=0.5$ and $p_1=0$, then $F$ is exactly the famed Cantor function (see \cite[Proposition 4.2]{Dovgoshey2006Cantor}). The second graph in Figure \ref{fig:unfair_example} showcases this example. Thus, these cumulative distribution functions can be thought of as probabilistic generalizations of the Cantor function.
	\end{exmp}
	
	\captionsetup[subfigure]{labelformat=empty}
	\captionsetup[subfigure]{font=tiny} 
	
	In the context of an actual multi-state device, this result says that unless your device is absolutely perfect and unbiased, its distribution will always be fractal-like and hence extremely difficult to work with. In light of these challenging distributions, given an imperfect device, it is desirable to quantify how close it is to an ideal, fair reference (i.e., how close the associated distribution, $F(x)$, is to the uniform distribution $y=x$). The authors of\cite{Pfeffer2025Dice} chose to do so by computing the unfair distribution's absolute deviation from fair:
	\begin{equation} \label{eq:metric}
		\|x-F(x)\|_\infty := \sup\{ |x- F(x)| \ : \ x\in [0,1]\}.
	\end{equation}
	This was chosen as the comparison metric since, from a statistical viewpoint, this value is equal to the well-known Kolmogorov-Smirnov (KS) statistic. While explicit forms for $\|x-F(x)\|_\infty$ are intractable due to the fractal properties of $F$, they do admit  tractable bounds. With $C[0,1]$ denoting the space of continuous functions on $[0,1]$, 
	\begin{theorem}[Theorem 3 in \cite{Pfeffer2025Dice}] \label{theorem: uniformcomp}
		For $q\in \{2,3,\ldots\}$, let $T\colon C[0,1]\to C[0,1]$ be the operator defined by
		\begin{equation} \label{seqdef}
			(Tf)(x) =   \begin{cases} 
				p_0f(qx) &\text{if} \ 0\leq x \leq \textstyle \frac{1}{q}\\
				p_0+p_1f(qx-1) &\text{if} \ \textstyle \frac{1}{q} \leq x \leq \frac{2}{q}\\
				\ \ \ \ \ \ \ \ \vdots &\ \ \ \ \ \ \ \ \ \vdots\\
				(p_0+p_1+\ldots+p_{q-2}) + p_{q-1}f(qx-(q-1)) &\text{if} \ \textstyle \frac{q-1}{q} \leq x \leq 1\\
			\end{cases}.
		\end{equation}
		Given the sequence of functions $(f_n)_{n=0}^\infty$ defined by $f_0=x$ and $f_{n+1}= Tf_n$, the distribution function $F$ in Theorem \ref{maintheorem1}, and $p_{\text{max}} = \text{max}\{p_0, \ldots, p_{q-1}\}$, it follows that $f_n \to F$ and 
		\begin{equation} \label{eq:bounds}
			\left(\frac{1}{1+p_{\text{max}}}\right)\|x-f_1(x)\|_\infty \leq \|x-F(x)\|_\infty \leq \displaystyle\left( \frac{1}{1-p_{\text{max}}}\right)\|x-f_1(x)\|_\infty.
		\end{equation}
	\end{theorem}
	This is notable since it estimates the quantity in \eqref{eq:metric} via $\|x-f_1\|_\infty$, where $f_1$ is a (finite) piece-wise linear function and has no fractal-like components.
	
	In this paper, we continue this investigation into the analytic properties of $F$. We make clear here that some of our results assume that the multi-state device is `typical' in the sense that $p_i\neq 0,1$ for all $i$, which is the modality in which most devices behave.  
	That said, the set of device configurations that are `atypical' is essentially a set of measure zero, so there is no harm in disregarding them when needed.
	
	To start, we again note that there are many ways to quantify the distance between $F$ and the ideal, uniform distribution $y=x$. Theorem \ref{theorem: uniformcomp} provided upper and lower bounds on $F$'s absolute deviation from fair (i.e., the sup-norm). Another option is to instead consider the $L^1$-norm that computes the absolute area between $F$ and $y=x$:
	\[
	\|x-F(x)\|_1 := \int_0^1 |x-F(x)|\mathop{dx}.
	\]
	Just as with the sup-norm from Theorem \ref{theorem: uniformcomp}, an explicit form for $\|x-F(x)\|_1$ is intractable, however, they admit easily computable upper and lower bounds. We show the following:
	
	\vspace{0.5pc}
	
	\noindent \textbf{Theorem \ref{thm:L^1Distance}.}
	If $F$ is the distribution defined in Theorem~\ref{maintheorem1} and $f_1$ is the first iterate of the operator $T$ defined in Theorem \ref{theorem: uniformcomp}, then 
	\[ 
	\left(\frac{q}{q+1}\right)A \leq \|x-F(x)\|_1 \leq \left(\frac{q}{q-1}\right)A,
	\] 
	where $A=\|x-f_1(x)\|$.
	
	\vspace{0.5pc}
	
\begin{figure*}[t]
		\centering
		
		\begin{subfigure}[b]{0.33\textwidth}
			\centering
			\includegraphics[width=\linewidth]{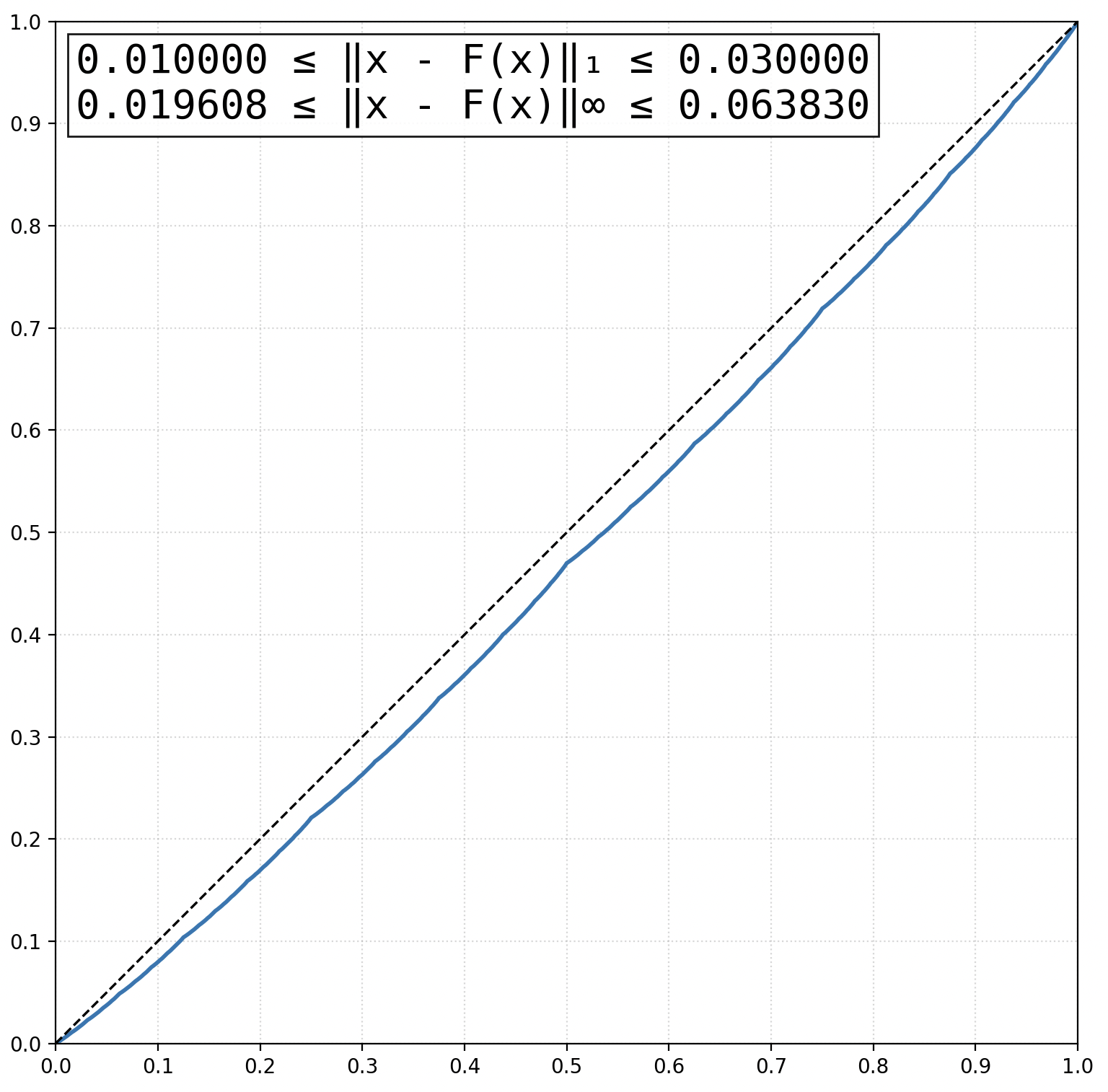}
			\caption{$P=\{0.47,\,0.53\}$}
			\label{fig:norm_1}
		\end{subfigure}\hfill
		\begin{subfigure}[b]{0.33\textwidth}
			\centering
			\includegraphics[width=\linewidth]{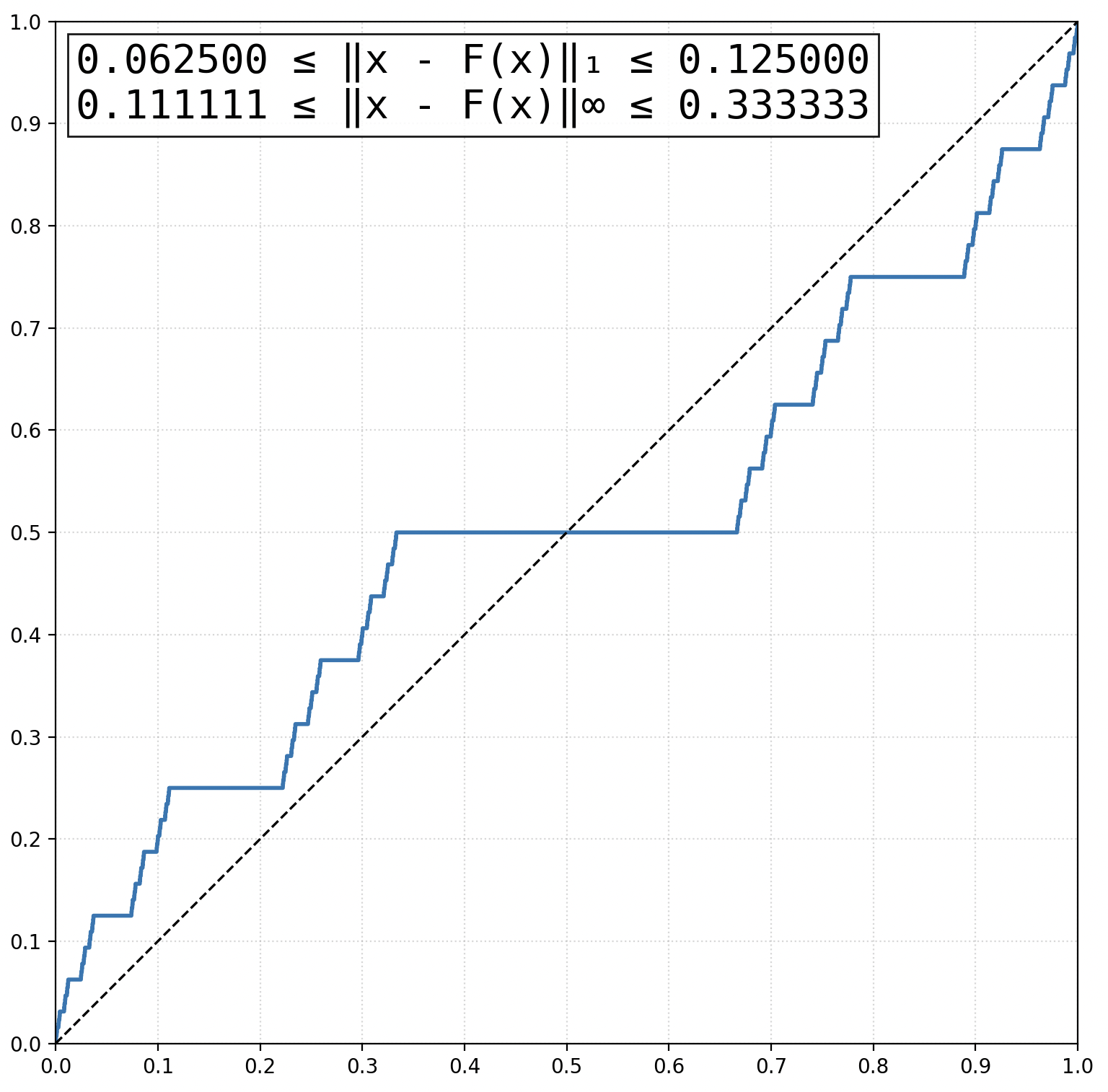}
			\caption{$P=\{0.50,\,0,\,0.50\}$}
			\label{fig:norm_2}
		\end{subfigure}\hfill
		\begin{subfigure}[b]{0.33\textwidth}
			\centering
			\includegraphics[width=\linewidth]{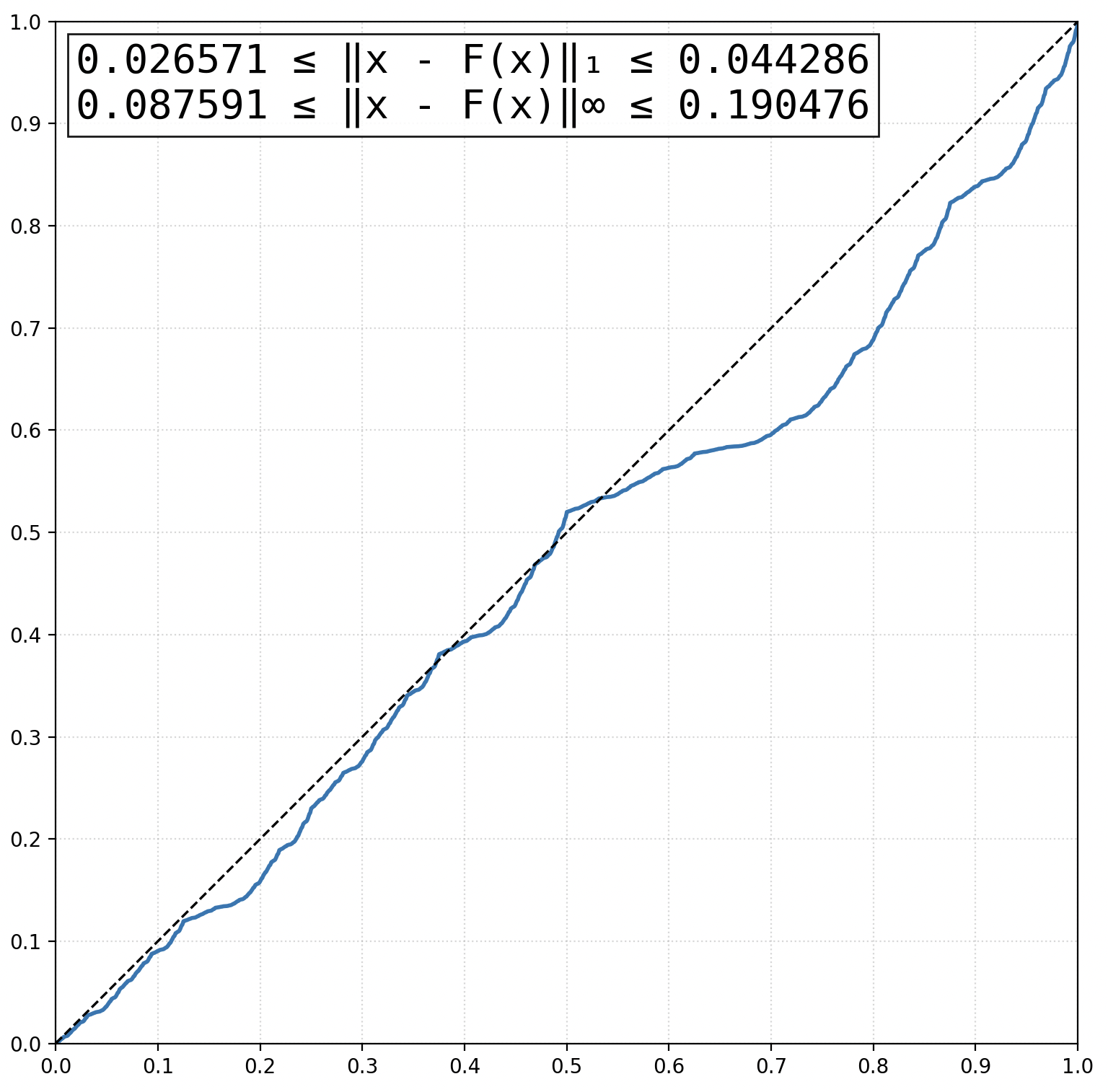}
			\caption{$P=\{0.23,\,0.29,\,0.11,\,0.37\}$}
			\label{fig:norm_3}
		\end{subfigure}\hfill
		\caption{The $L_1$ and sup-norm comparison for three different distributions, where the associated probability weights are given as $P=\{p_i\}$.}
		\label{fig:norm_comp}
	\end{figure*}

	As in Theorem \ref{theorem: uniformcomp}, $f_1$ is easy to work with since it is a (finite) piece-wise linear function. While the bounds for the sup-norm given in  \eqref{eq:bounds} are dependent on $p_{\text{max}}$ and blow up when the multi-state device is heavily skewed toward one state, the above bounds for the $L^1$-norm don't exhibit this behavior and, in fact, tighten rapidly with an increase in the number of sides on our die, $q$. See Figure \ref{fig:norm_comp} for examples of distributions and their associated norm bounds.

	Since the distributions $F$ are fractal-like, it is natural to ask questions about their smoothness. As noted in (i) of Theorem \ref{maintheorem1}, it is known that $F$ is continuous. We show that this can be strengthened to H\"{o}lder continuous:
	
	\vspace{0.5pc}
	
	\noindent \textbf{Theorem \ref{thm:holder1}.}
	Let $q\geq 2$, $p_i\neq 1$ for all $i$, and $p_{\text{max}} =\max\{p_i\}$. If $F$ is the distribution function defined in Theorem~\ref{maintheorem1}, then $F$ is H\"{o}lder continuous insofar that, given $x,y\in [0,1]$,
	\[
	|F(x)-F(y)| \leq (q+1)|x-y|^\alpha,
	\]
	where $\alpha= \frac{\ln(p_{\text{max}})}{  \ln\left(\frac{1}{q}\right)}$. Moreover, the exponent $\alpha$ is sharp.
	
	\vspace{0.5pc}
	
	\noindent \textbf{Theorem \ref{thm:holder2}.}
	Let $q\geq 2$, $p_i\neq 0$ for all $i$ and  $p_{\text{min}}=\min\{p_i\}$. If $F$ is the distribution function defined in Theorem~\ref{maintheorem1}, then, $F$ is lower-H\"{o}lder continuous insofar that, given $x,y\in [0,1]$,
	\[
	\left( \frac{1}{2q} \right)^\beta |x-y|^\beta \leq |F(x)-F(y)|,
	\]
	where $\beta= \frac{\ln(p_{\text{min}})}{  \ln\left(\frac{1}{q}\right)}$. Moreover, the exponent $\beta$ is sharp.
	
	\vspace{0.5pc}
	
	An immediate application of these results help provide upper and lower bounds on the Hausdorff dimension of the fractal-like support of $dF$:
	
	\vspace{0.5pc}
	
	\noindent \textbf{Theorem \ref{thm:hausdorffbounds}.}
	Let $F$ be the distribution defined in Theorem \ref{maintheorem1}, $p_i\neq0,1$ for all $i$, $p_{\text{min}}=\min\{p_i\}$, and $p_{\text{max}}=\max\{p_i\}$. If $\dim_\text{H}(\text{supp}(dF))$ denotes the Hausdorff dimension of $\text{supp}(dF)$, then
	\[
	\frac{\ln(p_{\text{max}}) }{ \ln\left(\frac{1}{q}\right)} \leq \dim_\text{H}(\text{supp}(dF)) \leq \frac{ \ln(p_{\text{min}}) }{\ln\left(\frac{1}{q}\right) }.
	\]
	
	\vspace{0.5pc}

	We now turn our attention toward approxiating $F$. If one were to stop the recursive definition for $F$ given in \eqref{qsideddierecursionmain} after $k$ iterations, the result would be a piecewise linear approximation to $F$ with $q^k$-many linear components. Moreover, this approximation can cross $F$ infinitely often. We provide a cheaper method of approximating $F$ by establishing a pair of step functions that bound $F$ from above and below, which converge uniformly to $F$ and only depend on computing the moments for $F$. Importantly, using our discovered H\"{o}lder exponents, we give an explicit stopping point for which these bounds are within a given error of $F$. While one can feasibly use empirical moments obtained directly from output data, we establish exact formulae for the moments of $F$ that depend only on $q$ and $p_j$. To this end, for $q\geq 2$ and $k\geq 0$, put $I_k(q) := \int_0^1 x^{k} F(x) \mathop{dx}$. We show the following:
	
	\vspace{0.5pc}
	
	\noindent \textbf{Theorem \ref{thm:l1}.}
	For $q\geq 2$, if $F$ is the distribution defined in Theorem \ref{maintheorem1}, then, for $k\geq 0$, the value of $I_k(q)$ is given recursively as
	\[
	I_k(q) = \int_0^1x^kF(x)\mathop{dx}=\frac{1}{q^{k+1}-1}\sum_{\ell=1}^{q-1}\left(\left(\frac{(\ell+1)^{k+1}-\ell^{k+1}}{k+1}\right)\cdot \sum_{j=0}^{\ell-1} p_j  + p_\ell\sum_{i=1}^k \binom{k}{i}\ell^i I_{k-i}(q)\right),
	\]
	where the last term is assumed to vanish when $k=0$.
	
	\vspace{0.5pc}
	
	\noindent \textbf{Corollary \ref{corr:k=rawmoment}.}
	For $q\geq 2$, if $F$ is the distribution defined in Theorem \ref{maintheorem1}, then, for $n\geq 1$, its $n$-th raw moment is given recursively by
	\[
	\mathbb{E}[X^n] = 1-\frac{n}{q^{n}-1}\sum_{\ell=1}^{q-1}\left(\left(\frac{(\ell+1)^{n}-\ell^{n}}{n}\right)\cdot \sum_{j=0}^{\ell-1} p_j  + p_\ell\sum_{i=1}^{n-1} \binom{n-1}{i}\ell^i I_{n-1-i}(q)\right),
	\]
	where the last term is assumed to vanish when $n=1$.
	
	\vspace{0.5pc}
	
	After the proof of Corollary \ref{corr:k=rawmoment}, we use these formulae to write down a recursive formula for the moment generating function $M_X(t)=\mathbb{E}[e^{tX}] = \sum_{n=0}^\infty \frac{t^n\mathbb{E}[X^n]}{n!}$. Additionally, as seen in Corollary \ref{corr:meanandvariance}, we also obtain closed formulae for the mean, $\mu(q)$, and variance, $\sigma^2(q)$, of $F$. 
	
\begin{figure}[h]
		\centering
		\begin{subfigure}[b]{0.33\textwidth}
			\centering
			\includegraphics[width=\linewidth]{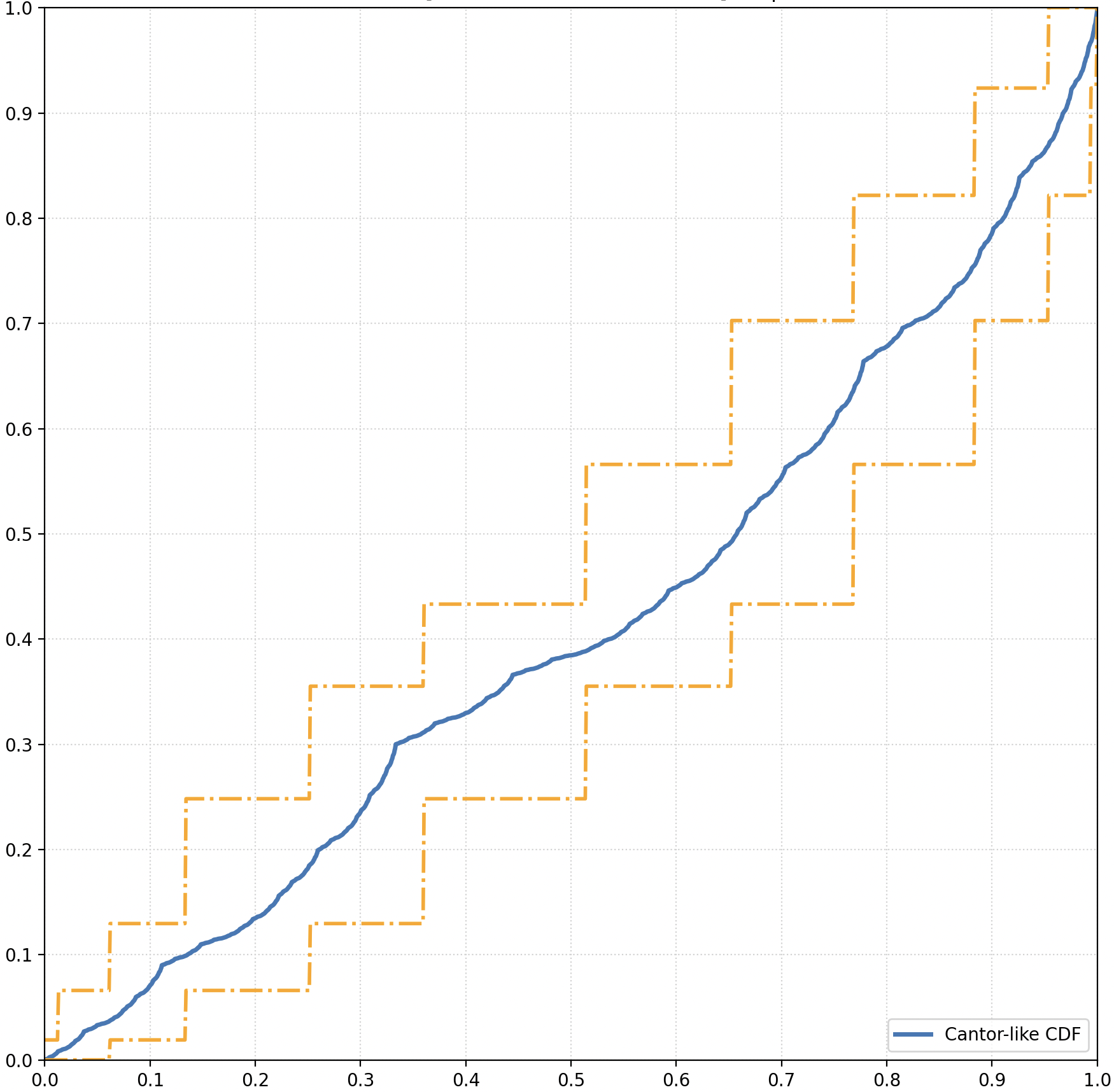}
			\caption{$P=\{0.48,\,0.22, 0.3\}$}
			\label{fig:hölder_1}
		\end{subfigure}\hfill
		\begin{subfigure}[b]{0.33\textwidth}
			\centering
			\includegraphics[width=\linewidth]{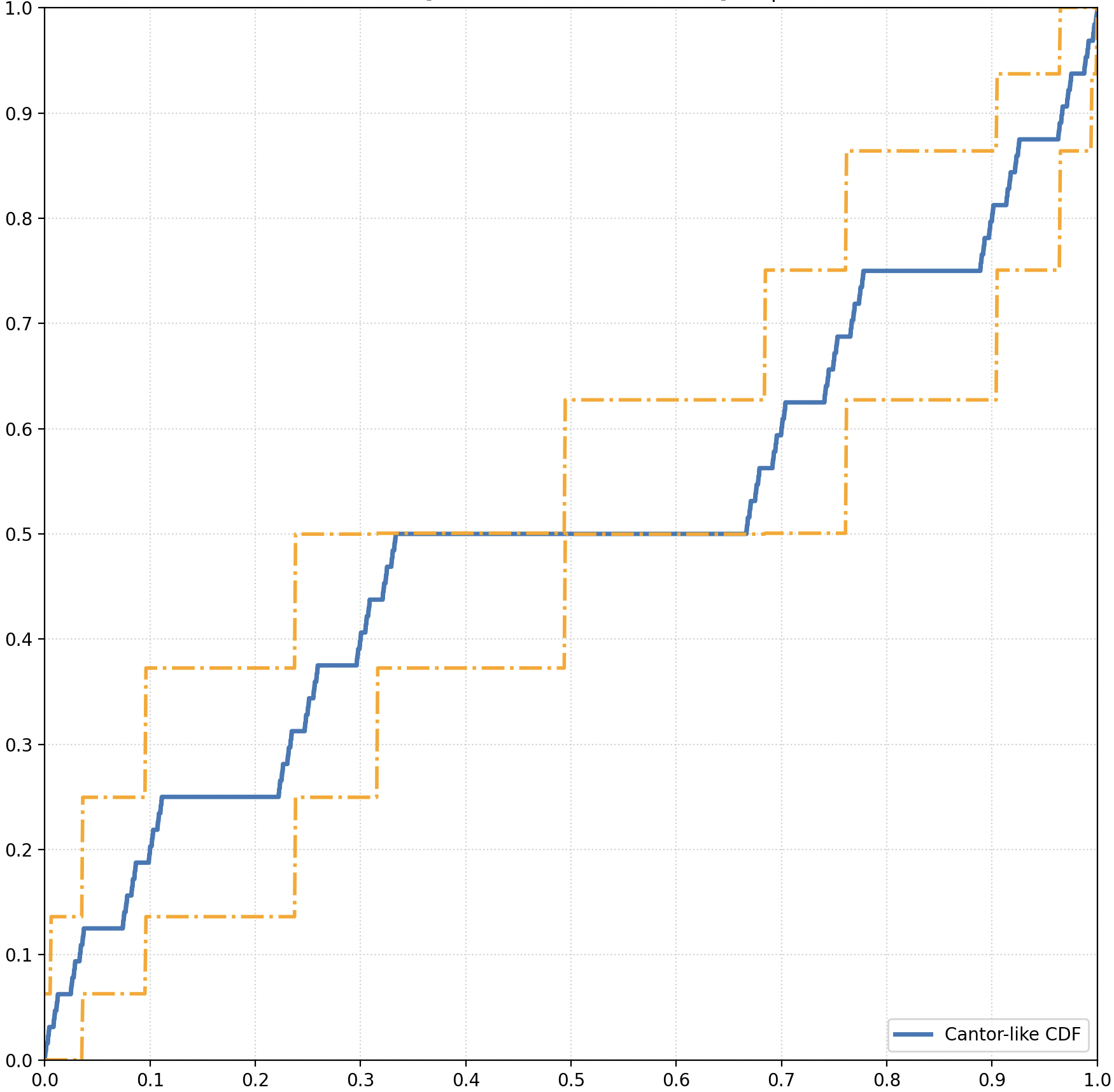}
			\caption{$P=\{1/2,\,0,\,1/2\}$}
			\label{fig:hölder_2}
		\end{subfigure}\hfill
		\begin{subfigure}[b]{0.33\textwidth}
			\centering
			\includegraphics[width=\linewidth]{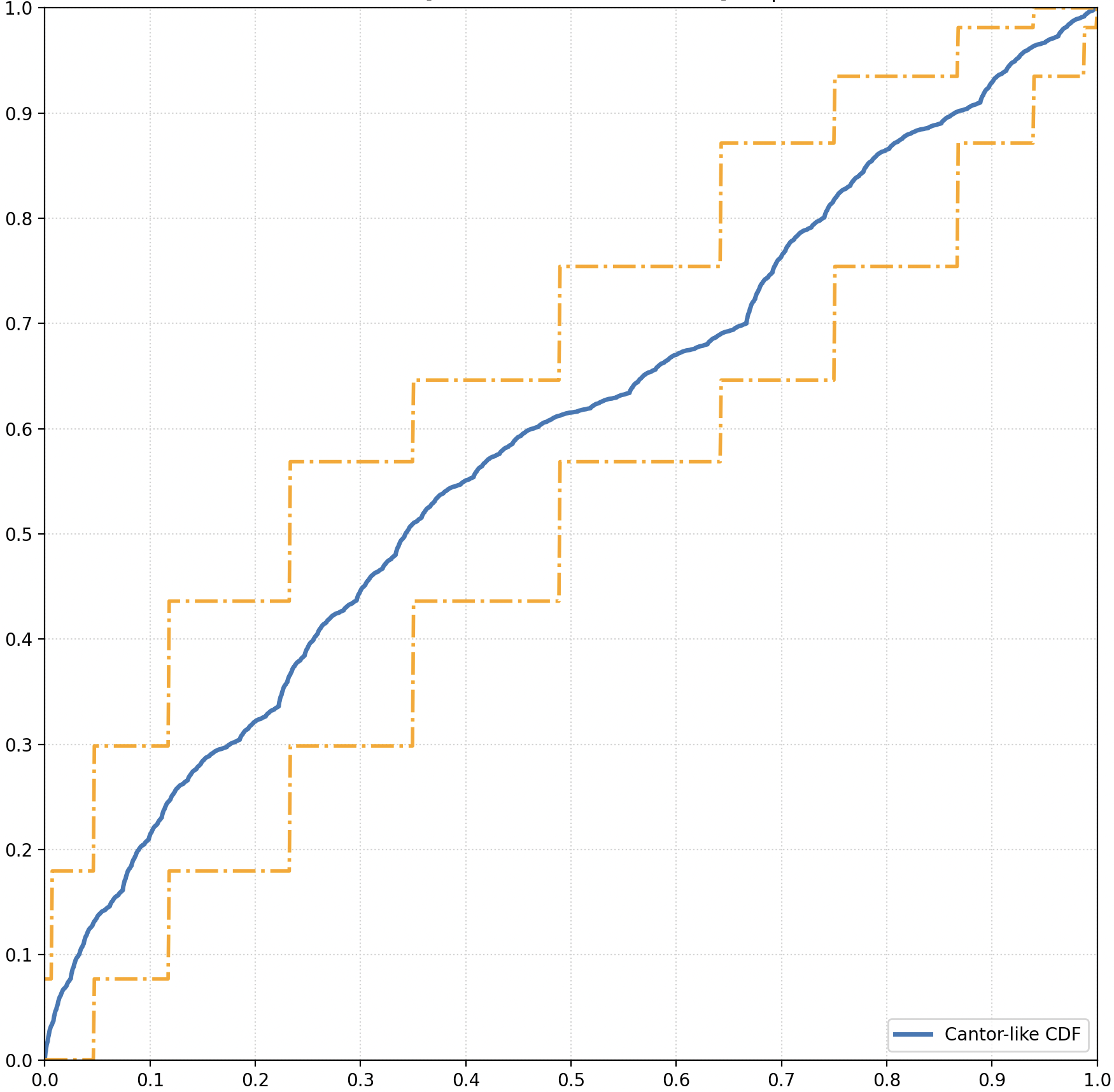}
			\caption{$P=\{0.3,\,0.22,\,0.48\}$}
			\label{fig:hölde_3}
		\end{subfigure}\hfill
		\caption{Examples of the bounds given by Theorem \ref{thm: moment bounds}. The distribution $F$ is given in blue, whereas the upper and lower bounds are given in orange. Both bounds are generated using $12$ moments. In each graph, the associated probability weights are given as $P=\{p_i\}$.}
		\label{fig:momentbounds}
	\end{figure}

	Finally, by appealing to the Chebyshev-Markov-Stieltjes inequalities (see Proposition \ref{prop:CMSineq}), the aforementioned step-function upper and lower approximations to $F$ are established (see Figure \ref{fig:momentbounds} for examples):
	
	\vspace{0.5pc}
	
	\noindent \textbf{Theorem \ref{thm: moment bounds}.}
	Let $F$ be the (continuous) cumulative distribution function of an arbitrary random variable $X$ defined on $[0,1]$ with moments $\{\mathbb{E}[X^k]\}_{k=0}^{2n-1}$ and let $\{P_i\}_{i=0}^n$ be the associated unique orthonormal polynomials.  If $\xi_1,\ldots,\xi_n$ are the distinct zeros of $P_n$, then
	\[
	\sum_{j=1}^{n-1} \lambda_n(\xi_j) \mathbbm{1}_{[\xi_{j+1},1]}(x) \leq F(x) \leq \lambda_n(\xi_1) + \sum_{j=1}^{n-1} \lambda_n(\xi_{j+1}) \mathbbm{1}_{(\xi_{j},1]}(x),
	\]
	where $\lambda_n(x)= \frac{1}{ \sum_{i=0}^{n-1} |P_i(x)|^2  }$ is the Christoffel function. Moreover, these bounds converge uniformly to $F$ from below and above, respectively, as $n\to \infty$ (i.e., as the number of moments considered increases).
	
	\vspace{0.5pc}

	Moreover, we provide the exact number of moments needed to approximate $F$ within a given error. To do so, Proposition \ref{prop:lambdatozero} establishes an explicit decay rate for $\lambda_n$ in terms of $\alpha= \frac{\ln(p_{\text{max}})}{  \ln\left(\frac{1}{q}\right)}$, the sharp H\"{o}lder exponent from Theorem \ref{thm:holder1}. Then, with
	\[
	L_n(x) = \sum_{j=1}^{n-1} \lambda_n(\xi_j) \mathbbm{1}_{[\xi_{j+1},1]}(x) \hspace{0.25in} \text{and} \hspace{0.25in} U_n(x) = \lambda_n(\xi_1) + \sum_{j=1}^{n-1} \lambda_n(\xi_{j+1}) \mathbbm{1}_{(\xi_{j},1]}(x)
	\]
	denoting the lower and upper bounds, respectively, in Theorem \ref{thm: moment bounds}, we obtain the following corollary:
	
	\vspace{0.5pc}
	
	\noindent \textbf{Corollary \ref{corr:momentstoboundF}.}
	Given $\varepsilon>0$, put $n=\left\lfloor \left(  \frac{8\cdot 16^\alpha (q+1)}{\varepsilon (4-\alpha)}\right)^{1/\alpha} \right\rfloor +2$, where $\alpha= \frac{\ln(p_{\text{max}})}{  \ln\left(\frac{1}{q}\right)}$ is the sharp H\"{o}lder exponent from Theorem \ref{thm:holder1}. If $L_n$ and $U_n$ are computed using the moments $\mathbb{E}[X^k]$ for $0\leq k\leq 2n-1$, then $\|U_n-L_n\|_\infty<\varepsilon$. In particular, so are $\|U_n-F(x)\|_\infty$ and $\|L_n-F(x)\|_\infty$.
	
	\vspace{0.5pc}
	
	\section{$L^1$ distance between $y=x$ and $F(x)$}
	
	\begin{theorem} \label{thm:L^1Distance}
		If $F$ is the distribution defined in Theorem~\ref{maintheorem1} and $f_1$ is the first iterate of the operator $T$ defined in Theorem \ref{theorem: uniformcomp}, then 
		%
		%
		\[ 
		\left(\frac{q}{q+1}\right)A \leq \|x-F(x)\|_1 \leq \left(\frac{q}{q-1}\right)A,
		\] 
		where $A=\|x-f_1(x)\|_1$.
	\end{theorem}
	\begin{proof}
		Consider $f_2=Tf_1$ and consider the intervals $\left[ 0, \textstyle \frac{1}{q} \right]$, $\left[ \frac{1}{q},\frac{2}{q}\right], \ldots, \left[ \frac{q-1}{q},1\right]$. Let $0\leq j\leq q-1$ be arbitrary and consider $\int_{j/q}^{(j+1)/q} |f_1(x)-f_2(x)|\mathop{dx}$. Applying the substitution $u=qx-j$, we get
		\begin{align*}
			\int_{j/q}^{(j+1)/q} |f_1(x)-f_2(x)|\mathop{dx} &= \int_{j/q}^{(j+1)/q} |Tx - Tf_1(x)|\mathop{dx}\\
			&= \int_{j/q}^{(j+1)/q} |(p_0+p_1+\ldots + p_j(qx-j) - (p_0+p_1+\ldots+p_{j}f_1(qx-j)|\mathop{dx}\\
			&= \int_{j/q}^{(j+1)/q} |p_j(qx-j) - p_{j}f_1(qx-j)|\mathop{dx}\\
			&= p_j\int_{j/q}^{(j+1)/q} |(qx-j) - f_1(qx-j)|\mathop{dx}\\
			&= \frac{p_j}{q}\int_{0}^{1} |u - f_1(u)|\mathop{du}\\
			&= \frac{p_j}{q}A\\
		\end{align*}
		Thus, summing across all intervals $\left[ 0, \textstyle \frac{1}{q} \right]$, $\left[ \frac{1}{q},\frac{2}{q}\right], \ldots, \left[ \frac{q-1}{q},1\right]$ yields
		\[
		\int_0^1 |f_1(x)-f_2(x)|\mathop{dx} = \sum_{j=0}^{q-1} \int_{j/q}^{(j+1)/q} |f_1(x)-f_2(x)|\mathop{dx}  = \sum_{j=0}^{q-1} \frac{p_j}{q}A = \frac{A}{q}(p_0+p_1+\ldots+p_{q-1}) = \frac{A}{q}.
		\]
		Iterating $T$ once more increases the number of intervals by a factor of $q$. In fact, similar reasoning as above shows that
		\begin{equation}\label{eq:iteratedTarea}
			\|Tf_n(x) - Tf_{n+1}(x)\|_1 = \frac{1}{q}\|f_n(x) - f_{n+1}(x)\|_1.
		\end{equation}
		Iterating multiple times yields
		\[
		A_n:=\int_0^1 |f_n(x) - f_{n+1}(x)| \mathop{dx} = \frac{A}{q^n},
		\]
		for $n\geq 1$ (where $A_0:=A$). Therefore, by Tonelli's theorem, 
		\[
		\int_0^1 |x-F(x)|\mathop{dx} \leq \sum_{n=0}^\infty A_n \leq \sum_{n=0}^\infty  \frac{A}{q^n}= A  \sum_{n=0}^\infty  \frac{1}{q^n} = \left(\frac{q}{q-1}\right)A.
		\]
		
		The lower bound is shown by invoking the reverse triangle inequality for norms. The following are the details: The triangle inequality for metric spaces guarantees that $\|f\|_1 = \|f-g+g\|_1 \leq \|f-g\|_1 + \|g\|_1$ for any $f,g\in C[0,1]$ so that $\|f\|_1-\|g\|_1 \leq \|f-g\|_1$. With $f_\infty := \lim f_n = F$, this shows that
		\begin{equation} \label{eq:lowerbound1}
			\| f_0-f_1 \|_1 - \|f_\infty - f_1\|_1 \leq \|(f_0-f_1)-(f_\infty - f_1)\|_1 = \|f_0-f_\infty\|_1.
		\end{equation}
		Since $f_\infty$ is a fixed point for $T$ (see the proof of Theorem \ref{theorem: uniformcomp} in \cite{Pfeffer2025Dice}), we know that $Tf_\infty = f_\infty$. Thus, \eqref{eq:iteratedTarea} shows that 
		\[
		\|f_\infty - f_1\|_1 = \|Tf_\infty-Tf_0\|_1 =\frac{1}{q}\|f_\infty-f_0\|_1.
		\]
		Negating both sides, we have $-\|f_\infty - f_1\|_1 = -\frac{1}{q}\|f_\infty-f_0\|_1$ so that $\|f_0-f_1\|_1-\|f_\infty - f_1\|_1 = \|f_0-f_1\|_1-\frac{1}{q}\|f_\infty-f_0\|_1$. Combining this with (\ref{eq:lowerbound1}), we have
		\[
		\|f_0-f_1\|_1-\frac{1}{q}\|f_\infty-f_0\|_1 = \|f_0-f_1\|_1-\|f_\infty - f_1\|_1 \leq \|f_0 - f_\infty\|_1.
		\]
		Solving the above for the quantity $\|f_0-f_\infty\|_1 = \|f_\infty-f_0\|_1$ yields
		\[
		\left( \frac{q}{1+q}\right) \|f_0-f_1\|_1 \leq \|f_0-f_\infty\|_1.
		\]
	\end{proof}
	
	\section{H\"{o}lder continuity for $F$}
	
	\begin{theorem} \label{thm:holder1}
		Let $q\geq 2$, $p_i\neq 1$ for all $i$, and $p_{\text{max}} =\max\{p_i\}$. If $F$ is the distribution function defined in Theorem~\ref{maintheorem1}, then $F$ is H\"{o}lder continuous insofar that, given $x,y\in [0,1]$,
		\[
		|F(x)-F(y)| \leq (q+1)|x-y|^\alpha,
		\]
		where $\alpha= \frac{\ln(p_{\text{max}})}{  \ln\left(\frac{1}{q}\right)}$. Moreover, the exponent $\alpha$ is sharp.
	\end{theorem}
	\begin{proof}
		Let $x,y\in [0,1]$ and, without loss of generality, suppose $y>x$. By the Archimedean property of the real line, there exists $n\in \mathbb{N}$ such that $\frac{1}{q^n} \leq  y-x < \frac{1}{q^{n-1}}$. Moreover, there exists $k\in \mathbb{N}$ such that $\frac{k}{q^n} \leq x < \frac{k+1}{q^n}$. Together, this guarantees that $
		y\geq x+\frac{1}{q^n} \geq \frac{k}{q^n} + \frac{1}{q^n} = \frac{k+1}{q^n}$ and $y < x+\frac{1}{q^{n-1}} = x+\frac{q}{q^n} < \frac{k+(q+1)}{q^n}$. Visually, we have the following:
		
		\vspace{0.5pc}
		
		\begin{center}
			\begin{tikzpicture}[scale=1,>=stealth,
				every node/.style={align=center,scale=1}] 
				\draw[black,line width=1pt, <->] (1,0)--(7,0);
				\draw (2,.1)--(2,-.1) node[below]{$\frac{k}{q^n}$};
				\draw (3,.1)--(3,-.1) node[below]{$x$};
				\draw (4,.1)--(4,-.1) node[below]{$\frac{k+1}{q^n}$};
				\draw (5,.1)--(5,-.1) node[below]{$y$};
				\draw (6,.1)--(6,-.1) node[below]{$\frac{k+(q+1)}{q^n}$};
			\end{tikzpicture}
		\end{center}  
		
		\vspace{0.5pc}
		
		Let $\gamma\leq \alpha= \frac{\ln(p_{\text{max}})}{  \ln\left(\frac{1}{q}\right)}$. Solving for $p_{\text{max}}$ gives $p_{\text{max}} \leq \left( \frac{1}{q}\right)^\gamma$. In the spirit of the proof of Theorem \ref{maintheorem1}, for  $1\leq \ell \leq q+1$, let $u_i^\ell \in \{0,1,\ldots q-1\}$ be the outcomes that uniquely identify the interval $\left[ \frac{k+(\ell-1)}{q^n}, \frac{k+\ell}{q^n}\right]$ insofar that
		\[
		F\left(\frac{k+\ell}{q^n}\right) - F\left( \frac{k+(\ell-1)}{q^n}\right) = p_{u_1^\ell}\cdot p_{u^\ell_2}\cdot\ldots\cdot p_{u_n^\ell},
		\]
		where $p_{u^\ell_i}$ is the probability that the outcome $u^\ell_i$ is witnessed. Since $F$ is increasing, we have
		\begin{align*}
			|F(y)-F(x)| &\leq \textstyle F\left(\frac{k+(q+1)}{q^n}\right) - F\left( \frac{k}{q^n}\right)\\
			&= \textstyle F\left(\frac{k+1}{q^n}\right) - F\left( \frac{k}{q^n}\right) + F\left(\frac{k+2}{q^n}\right) - F\left( \frac{k+1}{q^n}\right) + \ldots + F\left(\frac{k+(q+1)}{q^n}\right) - F\left( \frac{k+q}{q^n}\right)\\
			& = (p_{u_1^1}\cdot p_{u^1_2}\cdot\ldots\cdot p_{u_n^1}) + (p_{u_1^2}\cdot p_{u^2_2}\cdot\ldots\cdot p_{u_n^2}) + \ldots + (p_{u_1^{q+1}}\cdot p_{u^{q+1}_2}\cdot\ldots\cdot p_{u_n^{q+1}})\\
			&\leq (q+1)p_{\text{max}}^n = (q+1)\left( \left( \frac{1}{q}\right)^\gamma\right)^n = (q+1) \left( \frac{1}{q^n}\right)^\gamma\\
			&\leq (q+1)|y-x|^\gamma,
		\end{align*}
		where we note that $F$ is naturally interpreted to be equal to 1 to the right of $[0,1]$. Therefore, $F$ is H\"{o}lder continuous with exponent $\gamma \leq \alpha=\frac{\ln(p_{\text{max}})}{  \ln\left(\frac{1}{q}\right)}$. 
		
		We claim that this exponent is sharp insofar that, if $\gamma>\alpha$, then the H\"{o}lder inequality doesn't hold. Specifically, we show that, given any $C>0$, we can find $x,y\in[0,1]$ so that $|F(x)-F(y)|>C|x-y|^\gamma$. To see this, suppose $\gamma>\alpha=\frac{\ln(p_{\text{max}})}{  \ln\left(\frac{1}{q}\right)}$. This implies that $p_{\text{max}}>\left( \frac{1}{q}\right)^\gamma$. In particular, $p_{\text{max}} q^\gamma>1$.

		Now, let $C>0$ be arbitrary. Choose $n\in \mathbb{N}$ large enough so that $\left( p_{\text{max}} q^\gamma \right)^n > C$. Let $k\in \mathbb{N}$ be chosen so that, over the interval $[x,y]:=\left[\frac{k}{q^{n}}, \frac{k+1}{q^{n}} \right]$,  
		\[
		|F(y)-F(x)| = F(y)-F(x) = F\left(\frac{k+1}{q^{n}}\right) - F\left( \frac{k}{q^{n}}\right) = p_{\text{max}}\cdot p_{\text{max}}\cdot\ldots\cdot p_{\text{max}} = p_{\text{max}}^{n},
		\]
		where the initial absolute value can be dropped since $F$ is increasing. It now follows that 
		\[
		|F(y)-F(x)| = p_{\text{max}}^{{n}} = \left( p_{\text{max}} q^\gamma \right)^n \left( \frac{1}{q^n} \right)^\gamma > C|x-y|^\gamma.
		\]
		Thus, for any $C$, we can find $x,y\in [0,1]$ so that $|F(y)-F(x)|>C|x-y|^\gamma$. 
	\end{proof}

	H\"{o}lder continuity is, in some sense, a measure of how rugged a function's graph looks. It is known that when the H\"{o}lder exponent, $\alpha$, is equal to 1, then the associated function is Lipschitz continuous and is hence smooth and well-behaved, whereas when $\alpha\in(0,1)$ gets closer to zero, the graph gets more and more rugged. This aligns nicely with our setting. According to Theorem \ref{thm:holder1}, our distribution's H\"{o}lder exponent is equal to $\alpha=  \ln(p_{\text{max}})/ \ln\left(\textstyle 1/q\right)$. This is a ratio between the worst probability weighting and the ideal, uniform weighting. When these values are equal (i.e., $p_{\text{max}}=\textstyle \frac{1}{q}$), this means $p_i=\frac{1}{q}$ for all $i$. This makes $F$ exactly the uniform distribution $y=x$, which exhibits no ruggedness. As the probability weights deviate from the fair $1/q$, the fractal-like nature of $F$ worsens, producing a more rugged graph and a smaller H\"{o}lder exponent.
	
	We pause to note that the constant $C=q+1$ in Theorem \ref{thm:holder1} is not sharp. The choices $q=3$, $p_0=p_2=0.5$, and $p_1=0$ all satisfy the hypotheses of Theorem \ref{thm:holder1} (so, $p_{\text{max}} = 0.5$) and, as noted in Example \ref{ex:cantor}, the distribution function $F$ in this case is exactly the Cantor function. In \cite[Proposition 10.1]{Dovgoshey2006Cantor}, it is shown that, while the exponent matches our result as $\alpha = \frac{\ln(1/2)}{\ln(1/3)}$, the sharp constant is equal to $C=1$ and not $C=q+1=4$.
	
	The following result shows that we can establish a lower-H\"{o}lder inequality as well:
	
	\begin{theorem} \label{thm:holder2}
		Let $q\geq 2$, $p_i\neq 0$ for all $i$ and  $p_{\text{min}}=\min\{p_i\}$. If $F$ is the distribution function defined in Theorem~\ref{maintheorem1}, then, given $x,y\in [0,1]$,
		\[
		\left( \frac{1}{2q} \right)^\beta |x-y|^\beta \leq |F(x)-F(y)|,
		\]
		where $\beta= \frac{\ln(p_{\text{min}})}{  \ln\left(\frac{1}{q}\right)}$. Moreover, the exponent $\beta$ is sharp.
	\end{theorem}
	\begin{proof}
		Let $x,y\in [0,1]$ and, without loss of generality, suppose $y>x$. Suppose further that $0<y-x<\frac{2}{q}$. By the Archimedean property of the real line, there exists $n\in \mathbb{N}$ such that $\frac{2}{q^n} \leq y-x < \frac{2}{q^{n-1}}$. Moreover, there exists $k\in \mathbb{N}$ such that $\frac{k}{q^n} \leq x < \frac{k+1}{q^n}$. Together, this guarantees that $
		y\geq x+\frac{2}{q^n} \geq \frac{k}{q^n} + \frac{2}{q^n} = \frac{k+2}{q^n}$. In particular, $[x,y]$ contains a base-$q$ interval. Visually, we have the following:
		
		\vspace{0.5pc}
		
		\begin{center}
			\begin{tikzpicture}[scale=1,>=stealth,
				every node/.style={align=center,scale=1}] 
				\draw[black,line width=1pt, <->] (0,0)--(6,0);
				\draw (1,.1)--(1,-.1) node[below]{$\frac{k}{q^n}$};
				\draw (2,.1)--(2,-.1) node[below]{$x$};
				\draw (3,.1)--(3,-.1) node[below]{$\frac{k+1}{q^n}$};
				\draw (4,.1)--(4,-.1) node[below]{$\frac{k+2}{q^n}$};
				\draw (5,.1)--(5,-.1) node[below]{$y$};
			\end{tikzpicture}
		\end{center}  
		
		\vspace{0.5pc}
		
		Again, let $u_i \in \{0,1,\ldots q-1\}$ be the outcomes that uniquely identify the interval $\left[ \frac{k+1}{q^n}, \frac{k+2}{q^n}\right]$ insofar that
		\[
		F\left(\frac{k+2}{q^n}\right) - F\left( \frac{k+1}{q^n}\right) = p_{u_1}\cdot p_{u_2}\cdot\ldots\cdot p_{u_n},
		\]
		where $p_{u_i}$ is the probability that the outcome $u_i$ is witnessed. It follows that
		\begin{align*}
			|F(y)-F(x)| &= F(y)-F(x)\geq F\left(\frac{k+2}{q^n}\right) - F\left( \frac{k+1}{q^n}\right)= p_{u_1}\cdot p_{u_2}\cdot\ldots\cdot p_{u_n}\\
			&\geq (p_{\text{min}})^n \geq \left(\frac{1}{q^n}\right)^\beta = \left(\frac{1}{2q}\right)^\beta \left( \frac{2q}{q^n}\right)^\beta \geq \left(\frac{1}{2q}\right)^\beta |x-y|^\beta,
		\end{align*}
		as desired. If, instead, $y-x\geq \frac{2}{q}$, then $[x,y]$ automatically contains a base-$q$ interval and the result follows via similar argumentation. Finally, the conclusion is trivial when $x=y$.
		
		We claim that this exponent is sharp insofar that, if $\gamma<\beta$, then, given any $C>0$, we can find $x,y\in[0,1]$ so that $|F(x)-F(y)|<C|x-y|^\gamma$. To see this, suppose $\gamma<\beta=\frac{\ln(p_{\text{min}})}{  \ln\left(\frac{1}{q}\right)}$. This implies that $p_{\text{min}}<\left( \frac{1}{q}\right)^\gamma$. Now, let $C>0$ be arbitrary. Choose $n\in \mathbb{N}$ large enough so that $\left( \frac{p_{\text{min}}}{(1/q)^\gamma} \right)^n < C$. Let $k\in \mathbb{N}$ be chosen so that, over the interval $[x,y]:=\left[\frac{k}{q^{n}}, \frac{k+1}{q^{n}} \right]$,  
		\[
		|F(y)-F(x)| = F(y)-F(x) = F\left(\frac{k+1}{q^{n}}\right) - F\left( \frac{k}{q^{n}}\right) = p_{\text{min}}\cdot p_{\text{min}}\cdot\ldots\cdot p_{\text{min}} = p_{\text{min}}^{{n}},
		\]
		where the initial absolute value can be dropped since $F$ is increasing. It now follows that 
		\[
		|F(y)-F(x)| = p_{\text{min}}^{{n}} = \left( \frac{p_{\text{min}}}{(1/q)^\gamma} \right)^n \left( \frac{1}{q^n} \right)^\gamma < C|x-y|^\gamma.
		\]
		Therefore, $|F(x)-F(y)|<C|x-y|^\gamma$. 
	\end{proof}
	
	These results allow us to discuss another topic of interest for fractal-like sets: the Hausdorff dimension of the support of $dF$. Theorem \ref{maintheorem1} observed that $F$ may be a singular distribution, meaning the measure $dF$ is supported on a set of Lebesgue measure zero. For example, in the case that $F$ is the Cantor function (see Example \ref{ex:cantor}), this support is known to be the Cantor set with Hausdorff dimension $\log_3(2)$, which is exactly equal to the Cantor function's H\"{o}lder exponent. While our results won't be able to provide such an exact value for the Hausdorff dimension of the support of an arbitrary $F$, we can provide upper and lower bounds on it. To proceed, we need the following lemmas:
	\begin{lemma}[Proposition 2.3 in \cite{falconer2014fractal} and Lemma 5.1.3 in \cite{BishopPeres2017}]\label{lemma:hausdorffupper}
		Let $X\subset \mathbb{R}^n$ and suppose that $f\colon X\to \mathbb{R}^m$ satisfies the H\"{o}lder condition that, given any $x,y\in X$, there exists a constant $C$ such that $|f(x)-f(y)|\leq c|x-y|^\alpha$. Then $\dim_\text{H}f(X) \leq (1/\alpha)\dim_{\text{H}}X$.
	\end{lemma}
	\begin{lemma} \label{lemma:hausdorfflower}
		Let $X\subset \mathbb{R}^n$ and suppose that $f\colon X\to \mathbb{R}^m$ is an injective function that satisfies the lower-H\"{o}lder condition that, given any $x,y\in X$, there exists a constant $C$ such that $|f(x)-f(y)|\geq c|x-y|^\beta$. Then $\dim_\text{H}f(X) \geq (1/\beta)\dim_{\text{H}}X$.
	\end{lemma}
	
	Note that the proof of Lemma \ref{lemma:hausdorfflower} is essentially the proof of Lemma \ref{lemma:hausdorffupper} but with the inequalities switched around (the injectivity hypothesis is necessary simply to make the lower-H\"{o}lder condition make sense). The following is our result:
	
	\begin{theorem} \label{thm:hausdorffbounds}
		Let $F$ be the distribution defined in Theorem \ref{maintheorem1}, $p_i\neq0,1$ for all $i$, $p_{\text{min}}=\min\{p_i\}$, and $p_{\text{max}}=\max\{p_i\}$. If $\dim_\text{H}(\text{supp}(dF))$ denotes the Hausdorff dimension of $\text{supp}(dF)$, then
		\[
		\frac{\ln(p_{\text{max}}) }{ \ln\left(\frac{1}{q}\right)} \leq \dim_\text{H}(\text{supp}(dF)) \leq \frac{ \ln(p_{\text{min}}) }{\ln\left(\frac{1}{q}\right) }.
		\]
	\end{theorem}
	\begin{proof}
		Choose $X=\text{supp}(dF)$ and $f=F$ so that $f(X)=[0,1]$, which has Hausdorff dimension equal to 1. Since $p_i\neq 0$ for all $i$, it follows from \eqref{item2thm1} of Theorem \ref{maintheorem1} that $F$ is injective. Let $\alpha$ and $\beta$ be those values from Theorems \ref{thm:holder1} and \ref{thm:holder2}, respectively. The desired result then follows from Lemmas \ref{lemma:hausdorffupper} and \ref{lemma:hausdorfflower}.
	\end{proof}
	
	Note that this result cannot recover the fact that the Hausdorff dimension for the Cantor set is equal to the H\"{o}lder exponent of the Cantor function since, unfortunately, the requirement that $p_1=0$ violates the hypotheses of Theorem \ref{thm:hausdorffbounds}. In fact, it is not even reasonable to make $p_1$ arbitrarily small to try and approximate the Cantor function since $\beta$ will shrink to 0 while $\alpha$ approaches the desired value of $\log_3(2)$ leading to the unhelpful inequality $0\leq \dim_\text{H} \leq \log_3(2)$. A separate argument using coverings would be necessary to achieve the desired lower bound.

	\section{Moments and step-function approximations to $F$}
	
	The moments of a probability distribution provide important information on their behavior. In this section we derive recursive formulae for computing the moments for the distributions $F$ defined in Theorem \ref{maintheorem1}. 
	
	For $n\geq 1$, the $n$-th moment of the cumulative distribution function $F$ is defined by the Riemann-Stieltjes integral $\mathbb{E}[X^n] = \text{R-S} \int_{-\infty}^\infty x^n \mathop{dF(x)}$, where $X$ is a random variable and $\mathbb{E}$ is the expectation operator. The moment generating function is then given by
	\[
	M_X(t)=\mathbb{E}[e^{tX}] = 1 + t\mathbb{E}[X] + \frac{t^2\mathbb{E}[X^2]}{2!}+\ldots+\frac{t^n\mathbb{E}[X^n]}{n!}+\ldots.
	\]
	In our setting, it will suffice to compute $\mathbb{E}[X^n] = \text{R-S}\int_0^1 x^n \mathop{dF(x)}$. We now aim to convert these Riemann-Stieltjes integrals into something more computable. To do so, we need to recall two known results:
	
	\begin{proposition}[Theorem 4.4.1 in \cite{burk2007garden}] \label{prop:burk1}
		If $f$ is continuous on $[a,b]$ and $\phi$ is monotone increasing on $[a,b]$, then 
		\[
		\text{R-S}\int_a^b f(x)\mathop{d\phi(x)}
		\]
		exists.
	\end{proposition}
	
	\begin{proposition}[Theorem 4.3.2 in \cite{burk2007garden}]\label{prop:burk2}
		Suppose $f$ and $\phi$ are bounded functions with no common discontinuities on the interval $[a,b]$, and that the Riemann-Stieltjes integral of $f$ with respect to $\phi$ exists. Then the Riemann-Stieltjes integral of $\phi$ with respect to $f$ exists, and
		\[
		\text{R-S}\int_a^b \phi(x)\mathop{df(x)} = f(b)\phi(b) - f(a)\phi(a) - \text{R-S}\int_a^b f(x)\mathop{d\phi(x)}.
		\]
	\end{proposition}
	
	It was seen in (i) of Theorem \ref{maintheorem1} that $f(x)=F(x)$ is continuous on $[0,1]$ (and hence, bounded). Moreover, $\phi(x)=x^n$ is evidently monotonically increasing (and bounded) on $[0,1]$ for $n\geq 1$. Thus, by Proposition \ref{prop:burk1}, $\text{R-S}\int_0^1 F(x)\mathop{dx^n}$ exists for all $n\geq 1$. Proposition \ref{prop:burk2} then guarantees that $\text{R-S} \int_0^1 x^n \mathop{dF(x)}$ exists for all $n\geq 1$ and
	\[
	\text{R-S}\int_0^1 F(x)\mathop{dx^n} = 1^nF(1) - 0^nF(0) - \text{R-S}\int_0^1 x^n\mathop{dF(x)} = 1- \text{R-S}\int_0^1 x^n\mathop{dF(x)}.
	\]
	Equivalently,
	\begin{equation} \label{eq:rstranslator}
		\text{R-S}\int_0^1 x^n\mathop{dF(x)} = 1-\text{R-S}\int_0^1 F(x)\mathop{dx^n} = 1- n\int_0^1 x^{n-1}F(x)\mathop{dx}.
	\end{equation}
	Therefore, in order to compute the moments for $F$, it suffices to understand the values of $\int_0^1 x^{n-1}F(x)\mathop{dx}$ for $n\geq 1$. To this end, for $q\geq 2$, put $I_k(q) := \int_0^1 x^{k} F(x) \mathop{dx}$.

	\begin{theorem}\label{thm:l1}
		For $q\geq 2$, if $F$ is the distribution defined in Theorem \ref{maintheorem1}, then, for $k\geq 0$, the value of $I_k(q)$ is given recursively as
		\[
		I_k(q) = \int_0^1x^kF(x)\mathop{dx}=\frac{1}{q^{k+1}-1}\sum_{\ell=1}^{q-1}\left(\left(\frac{(\ell+1)^{k+1}-\ell^{k+1}}{k+1}\right)\cdot \sum_{j=0}^{\ell-1} p_j  + p_\ell\sum_{i=1}^k \binom{k}{i}\ell^i I_{k-i}(q)\right),
		\]
		where the last term is assumed to vanish when $k=0$.
	\end{theorem}
	\begin{proof}
		As noted in \eqref{qsideddierecursionmain}, $F$ is given by 
		\[
		F(x) = 
		\begin{cases} 
			p_0F(qx) &\text{if} \ 0\leq x \leq \textstyle \frac{1}{q}\\
			p_0+p_1F(qx-1) &\text{if} \ \textstyle \frac{1}{q} \leq x \leq \frac{2}{q}\\
			\ \ \ \ \ \ \ \ \vdots &\ \ \ \ \ \ \ \ \ \vdots\\
			(p_0+p_1+\ldots+p_{q-2}) + p_{q-1}F(qx-(q-1)) &\text{if} \ \textstyle \frac{q-1}{q} \leq x \leq 1\\
		\end{cases}.
		\]
		Thus, we split the desired integral across the base-$q$ intervals $\left[0,\textstyle \frac{1}{q}\right]$, $\left[\textstyle \frac{1}{q},\frac{2}{q}\right]$, $\ldots$, $\left[ \textstyle \frac{q-1}{q}, 1\right]$. To this end, for $0\leq \ell \leq q-1$, we consider the $\ell$-th such interval $\left[ \textstyle \frac{\ell}{q},\frac{\ell+1}{q}\right]$ and apply the substitution $u=qx-\ell$ to yield
		\begin{align*}
			\int_{\frac{\ell}{q}}^{\frac{\ell+1}{q}}
			x^kF(x)\mathop{dx} &= \int_{\frac{\ell}{q}}^{\frac{\ell+1}{q}}
			x^k\left(
			(p_0+\ldots+p_{\ell-1}) + p_{\ell}F(qx-\ell)
			\right)\mathop{dx}\\
			&= (p_0+\ldots+p_{\ell-1})\cdot \frac{(\ell+1)^{k+1}-\ell^{k+1}}{(k+1)\cdot q^{k+1}} + \frac{p_{\ell}}{q^{k+1}}\int_0^1 (u+\ell)^kF(u)\mathop{du}\\
			&= \left(\frac{(\ell+1)^{k+1}-\ell^{k+1}}{(k+1)\cdot q^{k+1}}\right)\cdot \sum_{j=0}^{\ell-1} p_j + \frac{p_{\ell}}{q^{k+1}}\int_0^1 (u+\ell)^kF(u)\mathop{du}, \label{eq:1}
		\end{align*}
		where we allow the natural observation that $ \sum_{j=0}^{\ell-1} p_j=0$ when $\ell=0$. We now unpack the last integral above. The Binomial Theorem guarantees that $(u+\ell)^k = \sum_{i=0}^k \binom{k}{i}u^{k-i}\ell^i$ so that
		\begin{align*}
			\int_0^1 (u+\ell)^kF(u)\mathop{du} 
			&= \int_0^1 \left(\sum_{i=0}^k \binom{k}{i}u^{k-i}\ell^i F(u)\right)\mathop{du}\\
			&= \sum_{i=0}^k \left(\binom{k}{i}\ell^i \int_0^1 u^{k-i}F(u)\mathop{du}\right)\\
			&=  \sum_{i=0}^k \binom{k}{i}\ell^i I_{k-i}(q)\\
			&=   I_k(q) +\sum_{i=1}^k \binom{k}{i}\ell^i I_{k-i}(q).
		\end{align*}
		Thus, for each $0\leq \ell \leq q-1$, we have
		\[
		\int_{\frac{\ell}{q}}^{\frac{\ell+1}{q}} x^kF(x)\mathop{dx} = \left(\frac{(\ell+1)^{k+1}-\ell^{k+1}}{(k+1)\cdot q^{k+1}}\right)\cdot \sum_{j=0}^{\ell-1} p_j  + \frac{p_\ell}{q^{k+1}}\left(\sum_{i=1}^k \binom{k}{i}\ell^i I_{k-i}(q)\right) + \frac{p_\ell}{q^{k+1}}I_k(q).
		\]
		Therefore, by summing across all such base-$q$ intervals, we have
		\begin{align*}
			I_k(q) &= \sum_{\ell=0}^{q-1}\int_{\frac{\ell}{q}}^{\frac{\ell+1}{q}} x^kF(x)\mathop{dx}\\
			&= \sum_{\ell=0}^{q-1}\left(\left(\frac{(\ell+1)^{k+1}-\ell^{k+1}}{(k+1)\cdot q^{k+1}}\right)\cdot \sum_{j=0}^{\ell-1} p_j  + \frac{p_\ell}{q^{k+1}}\left(\sum_{i=1}^k \binom{k}{i}\ell^i I_{k-i}(q)\right)\right) + \sum_{\ell=0}^{q-1} \frac{p_\ell}{q^{k+1}}I_k(q).
		\end{align*}
		Since $p_0+p_1+\ldots+p_{q-1}=1$, the last integral in the above expression can be simplified:
		\[
		\sum_{\ell=0}^{q-1} \frac{p_\ell}{q^{k+1}}I_k(q) = \frac{1}{q^{k+1}}I_k(q)\sum_{\ell=0}^{q-1} p_\ell = \frac{1}{q^{k+1}}I_k(q).
		\]
		Thus, by substituting this formulation into the previous expression, we can solve for $I_k(q)$ to see that
		\begin{align*}
			I_k(q) &= \left(\frac{q^{k+1}}{q^{k+1}-1}\right) \sum_{\ell=0}^{q-1}\left(\left(\frac{(\ell+1)^{k+1}-\ell^{k+1}}{(k+1)\cdot q^{k+1}}\right)\cdot \sum_{j=0}^{\ell-1} p_j  + \frac{p_\ell}{q^{k+1}}\sum_{i=1}^k \binom{k}{i}\ell^i I_{k-i}(q)\right)\\
			&= \frac{1}{q^{k+1}-1}\sum_{\ell=0}^{q-1}\left(\left(\frac{(\ell+1)^{k+1}-\ell^{k+1}}{k+1}\right)\cdot \sum_{j=0}^{\ell-1} p_j  + p_\ell\sum_{i=1}^k \binom{k}{i}\ell^i I_{k-i}(q)\right),
		\end{align*}
		where, again, we allow the natural conclusion that $ \sum_{j=0}^{\ell-1} p_j=0$ when $\ell=0$. In fact, with this consideration, we see that the $\ell=0$ term of the above sum is equal to 0 and therefore we have
		\[
		I_k(q) = \frac{1}{q^{k+1}-1}\sum_{\ell=1}^{q-1}\left(\left(\frac{(\ell+1)^{k+1}-\ell^{k+1}}{k+1}\right)\cdot \sum_{j=0}^{\ell-1} p_j  + p_\ell\sum_{i=1}^k \binom{k}{i}\ell^i I_{k-i}(q)\right),
		\]
		where the last term is assumed to vanish when $k=0$.
	\end{proof}
	
	\begin{corollary} \label{corr:k=rawmoment}
		For $q\geq 2$, if $F$ is the distribution defined in Theorem \ref{maintheorem1}, then, for $n\geq 1$, its $n$-th raw moment is given recursively by
		\[
		\mathbb{E}[X^n] = 1-\frac{n}{q^{n}-1}\sum_{\ell=1}^{q-1}\left(\left(\frac{(\ell+1)^{n}-\ell^{n}}{n}\right)\cdot \sum_{j=0}^{\ell-1} p_j  + p_\ell\sum_{i=1}^{n-1} \binom{n-1}{i}\ell^i I_{n-1-i}(q)\right),
		\]
		where the last term is assumed to vanish when $n=1$.
	\end{corollary}
	\begin{proof}
		By definition, $\mathbb{E}[X^n] = \text{R-S}\int_0^1 x^n \mathop{dF(x)}$. The result follows from combining the formula given in Theorem \ref{thm:l1} with the reduction given in \eqref{eq:rstranslator}.
	\end{proof}

	This allows us to generate a recursive formula for the moment generating function as well:
	\begin{align*}
		M_X(t)&=\mathbb{E}[e^{tX}] = 1 + t\mathbb{E}[X] + \frac{t^2\mathbb{E}[X^2]}{2!}+\ldots+\frac{t^n\mathbb{E}[X^n]}{n!}+\ldots\\
		&= 1+\sum_{n=1}^\infty\left( \frac{t^n}{n!}\left(1-\frac{n}{q^{n}-1}\sum_{\ell=1}^{q-1}\left(\left(\frac{(\ell+1)^{n}-\ell^{n}}{n}\right)\cdot \sum_{j=0}^{\ell-1} p_j  + p_\ell\sum_{i=1}^{n-1} \binom{n-1}{i}\ell^i I_{n-1-i}(q)\right)\right)\right).
	\end{align*}
	
	Of course, specific moments give valuable statistical information:
	\begin{corollary} \label{corr:meanandvariance}
		For $q\geq 2$, let $F$ be the distribution defined in Theorem \ref{maintheorem1}. If $\mu(q)$ and $\sigma^2(q)$ denote $F$'s mean and variance, respectively, then
		\begin{enumerate}[(i)]
			\item $ \displaystyle \mu(q) = 1- \frac{1}{q-1} \sum_{i=0}^{q-2} (q-1-i)p_i$
			\item $\displaystyle \sigma^2(q) = 1- 2\left(\frac{1}{q^2-1}\left( \frac{1}{2}\sum_{j=0}^{q-2}\left(q^2-(j+1)^2\right)p_j + I_1(q) \sum_{i=1}^{q-1} i p_i\right)\right) - \left(1- \frac{1}{q-1} \sum_{\ell=0}^{q-2} (q-1-\ell)p_\ell\right)^2$. 
		\end{enumerate}
	\end{corollary}
	\begin{proof}
		The result follows from Corollary \ref{corr:k=rawmoment} and the fact that  $\mu(q) = \mathbb{E}[X]$ and $\sigma^2(q) = \mathbb{E}[X^2] - \mathbb{E}[X]^2$. Note that the first summation for $\sigma^2(q)$ looks different from the formula in Theorem \ref{thm:l1} due to a telescoping phenomenon.
	\end{proof}

	We now aim to use the moments for $F$ and the H\"{o}lder exponents seen in Theorems \ref{thm:holder1} and \ref{thm:holder2} to provide an alternative method to approximate $F$ with step-functions instead of the piecewise linear approximant guaranteed by the recursive process given in \eqref{qsideddierecursionmain}. The main tool here is the Chebyshev-Markov-Stieltjes inequalities. We are aware that much of this is standard affair in the theory of orthogonal polynomials and their relation to the moment problem, but we include some details in order to keep the discussion as self-contained as possible.
	
	For $n\geq 0$, let  $m_k:=\mathbb{E}[X^k]$ for $0\leq k \leq 2n-1$. 
	We begin by determining the unique orthonormal polynomials, $P_0(x), P_1(x),\ldots, P_n(x)$, associated with the distribution $F$. That is, polynomials for which
	\[
	\int_0^1 P_i(x)P_j(x) \mathop{dF(x)} = \delta_{j,i},
	\]
	where $\delta_{j,i}$ is the Kronecker delta function. There are a few equivalent ways to do this. The first is to take the standard monomial basis $\{1,x,x^2,\ldots,x^n\}$ and apply the Gram-Schmidt Orthogonalization algorithm to it with respect to the weighted inner-product space $L^2(dF)$, followed by normalization of each polynomial. From this, for example, we see clearly that $P_0(x)=1$. Another route is to solve a three-term recurrence 
	(see \cite{ Liu2021RecurrenceCoefficients} for details). A third method is by way of Hankel matrices. Specifically, for $1\leq \ell \leq n$, we form the $\ell$-th Hankel matrix of the moments,
	\[
	H_\ell := \begin{bmatrix}
		m_0 & m_1 & \cdots & m_{\ell-1}\\
		m_1 & m_2 & \cdots & m_\ell\\
		\vdots & \vdots & \ddots & \vdots \\
		m_{\ell-1} & m_\ell & \cdots & m_{2\ell-2}
	\end{bmatrix}.
	\]
	and denote its determinant by $\Delta_\ell = \det(H_\ell)$. Theorem 2.1 in \cite{Gautschi2004Orthogonal} guarantees that, for $1\leq \ell \leq n$, if
	\[
	Q_\ell(x) := \frac{1}{\Delta_\ell}\det \begin{bmatrix}
		m_0 & m_1 & m_2 & \cdots & m_\ell \\
		m_1 & m_2 & m_3 & \cdots & m_{\ell+1}\\
		\vdots & \vdots & \vdots & \ddots & \vdots \\
		m_{\ell-1} & m_\ell & m_{\ell+1} & \cdots & m_{2\ell -1}\\
		1 & x & x^2 & \cdots & x^\ell
	\end{bmatrix},
	\]
	then $P_\ell(x) = \textstyle \frac{ Q_{\ell}(x)}{ \|Q_{\ell}(x)\|_{L^2(dF)}}$. With the orthonormal polynomials in hand, we let $\xi_1, \ldots, \xi_n$ be the $n$ real roots of $P_n(x)$. It is worth noting that the $\xi_j$ all lie in the interior of the convex hull of $\text{supp}(dF)$ (see \cite[pg.1]{kozhan2026zeros})---specifically, $\xi_j \in (0,1)$. Now consider the Christoffel-Darboux kernel (see \cite[pg.6]{NEVAI19863} or \cite{Swiderski2022Christoffel}):
	\[
	K_n(x,y) = \sum_{i=0}^{n-1} P_i(x)P_i(y).
	\]
	This is the reproducing kernel for the $\mathcal{P}_{n-1}(dF)$, the subspace of $L^2(dF)$ consisting of polynomials of degree at most $n-1$. The Christoffel function is then 
	\[
	\lambda_n(x):= \frac{1}{K_n(x,x)} = \frac{1}{ \sum_{i=0}^{n-1} |P_i(x)|^2  }.
	\]
	Note that, by setting $\rho(x)=1$ in \cite[Theorem 3.4.1]{Szego1939}, it follows that $\lambda_n(\xi_1)+\ldots+\lambda_n(\xi_{n})=1$. Additionally, since $P_0(x)=1$, we see that $\lambda_1(x)=1$ and $\lambda_n(x)$ is non-increasing with $n$. In fact, \cite[Theorem 9.7]{MR2500498} observes that, if $\mu$ is a measure on $\mathbb{R}$ of compact support, then, for all $x\in \mathbb{R}$, 
	\[
	\lim_{n\to \infty} \lambda_n(x) = \mu(\{x\}).
	\]
	In our setting, Theorem \ref{maintheorem1} guaranteed not only that  $dF$ has compact support, but also that  $F$ was both singular and continuous and hence $\lim_{n\to \infty} \lambda_n(x) = \mu(\{x\}) = 0$ for all $x\in [0,1]$. That is, the Christoffel function $\lambda_n$ converges pointwise to the zero function. We aim to establish an explicit upper bound on this decay rate. To do so, we need a lemma:
	
	\begin{lemma} \label{lemma:fejer}
		For $n\geq 4$ and $x_0\in [0,1]$, there exists a polynomial $q_n \in \mathcal{P}_n$ such that $q_n(x_0)=1$, $0\leq q_n(y) \leq 1$ for all $y\in [0,1]$, and $|q_n(y)| \leq \frac{4096}{n^4|y-x_0|^4}$ for $y\neq x_0$.
	\end{lemma}
	\begin{proof}
		
		For $m\geq 1$ and $\theta\notin 2\pi\mathbb{Z}$, let
		\[
		J_m(\theta) = \frac{1}{m^4}\left( \frac{ \sin(m\theta/2)}{\sin(\theta/2)}\right)^4
		\]
		be the (normalized) fourth-power Jackson kernel. As noted in \cite[p.~54]{lorentz1986approximation}, $J_m$ is actually an algebraic polynomial in $\cos(\theta)$ of degree $2m-2$. Consider the substitution $\cos(\theta)=1-\frac{1}{2}(y-x_0)^2$, which is well-defined since, when $y,x_0\in [0,1]$, we have $\theta\in [0,\textstyle \frac{\pi}{3}]$. For $n\geq 4$, set $m=\lfloor n/4 \rfloor +1$ and let $q_n(y)$ be $J_m(\theta)$ with this substitution applied. 
		
		First observe that, since $J_m$ is a polynomial of degree $2m-2$ in $\cos(\theta)$ and $\cos(\theta)=1-\frac{1}{2}(y-x_0)^2$ is degree 2 in $y$, we have that $q_n(y)$ is a polynomial of degree $4(m-1)=4\lfloor n/4 \rfloor \leq n$ and hence $q_n\in \mathcal{P}_n$. Second, observe that when $y=x_0$, we have $\cos(\theta)=1$ so that $\theta=0$. Thus, since $J_m(\theta)$ is a continuous function, $q_n(x_0) = J_m(0) = \lim_{\theta \to 0} J_m(\theta)=1$. Next, the elementary inequality $|\sin(m\theta/2)|\leq m|\sin(\theta/2)|$ guarantees that $0\leq J_m(\theta)\leq 1$ for all $\theta$. This implies that $0\leq q_n(y) \leq 1$ for all $y\in [0,1]$.
		
		Finally, applying the trigonometric identity $\sin^2\left( \textstyle \frac{\theta}{2}\right) = \frac{1-\cos(\theta)}{2}$ to the substitution $\cos(\theta)=1-\frac{1}{2}(y-x_0)^2$ gives $\sin\left( \textstyle \frac{\theta}{2}\right) = \frac{|y-x_0|}{2}$. Thus, since $|\sin(m\theta/2)|\leq 1$ and $m>\frac{n}{4}$, we have, for $y\neq x_0$,
		\[
		0\leq q_n(y) \leq \frac{1}{m^4 \left(\frac{ |y-x_0|}{2}\right)^4} = \frac{16}{m^4 |y-x_0|^4} < \frac{4096}{n^4|y-x_0|^4}.
		\]
	\end{proof}
	
	\begin{proposition} \label{prop:lambdatozero}
		If $q\geq 2$, $F$ is the distribution defined in Theorem \ref{maintheorem1}, and $\alpha$ is the sharp H\"{o}lder exponent for $F$ given in Theorem \ref{thm:holder1}, then, for any $x_0\in[0,1]$, it follows that $\lambda_1(x_0)=1$ and $\lambda_{n+1}(x_0) \leq \frac{4\cdot 16^\alpha (q+1)}{n^\alpha(4-\alpha)}$ for $n\geq 1$. In particular, $\lambda_{n+1}(x_0)\to 0$ as $n\to \infty$.
	\end{proposition}
	\begin{proof}
		Let $n\geq 0$. As noted earlier, $\lambda_{n+1}(x_0)$ is decreasing in $n$ and $\lambda_1(x_0)=1$. Therefore, when $1\leq n \leq 7$, it can be verified directly that $\lambda_{n+1}(x_0)\leq \lambda_1(x_0)=1 \leq \frac{4\cdot 16^\alpha (q+1)}{n^\alpha(4-\alpha)}$. Thus, we assume that $n\geq 8$. As noted in \cite[pg.9]{NEVAI19863} or \cite{Swiderski2022Christoffel}), we have the following alternative representation for $\lambda_{n+1}$:
		\[
		\lambda_{n+1}(x) = \inf_{p\in \mathcal{P}_{n}}\left\{ \int_0^1 |p(y)|^2\mathop{dF(y)}: p(x)=1\right\}.
		\]
		Let $x_0\in [0,1]$. Lemma \ref{lemma:fejer} guarantees that there exists a polynomial $q_n \in \mathcal{P}_n$ such that $q_n(x_0)=1$, $0\leq q_n(y)\leq 1$ for all $y\in [0,1]$, and $|q_n(y)| \leq \frac{4096}{n^4|y-x_0|^4}$ for $y\neq x_0$. In particular, this implies that $|q_n(y)|^2 = q_n(y)^2 \leq q_n(y) \leq \frac{4096}{n^4|y-x_0|^4}$ as well. Therefore,
		\begin{align*}
			\lambda_{n+1}(x_0) \leq \int_0^1 |q_n(y)|^2 \mathop{dF(y)} &= \int_{|y-x_0| \leq \textstyle \frac{8}{n}} |q_n(y)|^2 \mathop{dF(y)} + \int_{|y-x_0|> \textstyle \frac{8}{n}} |q_n(y)|^2 \mathop{dF(y)}\\
			&\leq \int_{|y-x_0| \leq \textstyle \frac{8}{n}} 1  \mathop{dF(y)} + \int_{|y-x_0|> \textstyle \frac{8}{n}} \frac{4096}{n^4|y-x_0|^4} \mathop{dF(y)}\\
			&=\int_0^{8/n} 1\mathop{d\nu(r)} + \int_{8/n}^1 \frac{4096}{n^4r^4} \mathop{d\nu(r)}\\
			&=\nu\left( \textstyle \frac{8}{n} \right)+ \int_{8/n}^1 \frac{4096}{n^4r^4} \mathop{d\nu(r)},
		\end{align*}
		where the penultimate equality comes from reframing the integrals to track the areas with respect to radius, $r$, from $x_0$ via $\nu(r):=dF([x_0-r,x_0+r])$ instead of tracking $y$.
		By Proposition \ref{prop:burk2} we have
		\begin{align*}
			\int_{8/n}^1 \frac{1}{n^4r^4} \mathop{d\nu(r)} &= \frac{\nu(1)}{n^4\cdot 1^4} - \frac{\nu(8/n)}{n^4(8/n)^4} + \int_{8/n}^1 \frac{4\nu(r)}{n^4 r^5}\mathop{dr}\\
			&=\frac{\nu(1)}{n^4} - \frac{\nu\left( \textstyle \frac{8}{n} \right)}{4096} + \int_{8/n}^1 \frac{4\nu(r)}{n^4 r^5}\mathop{dr}.
		\end{align*}
		Note that $\nu(1)$ makes sense if we make the natural assumption that $dF$ has no support beyond $[0,1]$. Recall, by Theorem \ref{thm:holder1},  $F$ is H\"{o}lder continuous with sharp exponent $\alpha<1$ and constant $q+1$ (again, the H\"{o}lder inequality still holds when $F$ is naturally extended beyond $[0,1]$). Therefore,
		\[
		\nu(r) = dF([x_0-r,x_0+r])=F(x_0+r)-F(x_0-r) \leq (q+1)((x_0+r)-(x_0-r))^\alpha = (q+1)2^\alpha r^\alpha.
		\]
		Hence, $\nu(1)\leq(q+1)2^\alpha$ and 
		\begin{align*}
			\int_{8/n}^1 \frac{4\nu(r)}{n^4 r^5}\mathop{dr} &\leq \frac{4\cdot 2^\alpha (q+1)}{n^4}\int_{8/n}^1 r^{\alpha-5}\mathop{dr}\\
			&=\frac{4\cdot 2^\alpha (q+1)}{n^4}\cdot \frac{1}{4-\alpha} \left(\frac{n^{4-\alpha}}{8^{4-\alpha}}-1\right)\\
			&= \frac{8^\alpha\cdot 4\cdot 2^\alpha (q+1)}{4096n^\alpha(4-\alpha)} - \frac{4\cdot 2^\alpha(q+1)}{(4-\alpha)n^4}.
		\end{align*}
		All together, we find that
		\begin{align*}
			\lambda_{n+1}(x_0) &\leq \nu\left( \textstyle \frac{8}{n} \right)+ \int_{8/n}^1 \frac{4096}{n^4r^4} \mathop{d\nu(r)}\\
			&\leq \nu\left( \textstyle \frac{8}{n} \right)+\frac{4096\nu(1)}{n^4} - \nu\left( \textstyle \frac{8}{n} \right) + \int_{8/n}^1 \frac{4096\cdot 4\nu(r)}{n^4 r^5}\mathop{dr}\\
			&\leq \frac{4096(q+1)2^\alpha}{n^4} + \frac{8^\alpha\cdot 4\cdot 2^\alpha (q+1)}{n^\alpha(4-\alpha)}- \frac{4096\cdot 4\cdot 2^\alpha(q+1)}{(4-\alpha)n^4}\\
			&= \frac{4\cdot 16^\alpha (q+1)}{n^\alpha(4-\alpha)} + \frac{4096(q+1)2^\alpha}{n^4}\left(1-\frac{4}{4-\alpha}\right)\\
			&\leq \frac{4\cdot 16^\alpha (q+1)}{n^\alpha(4-\alpha)},
		\end{align*}
		where the last inequality follows from the fact that $\alpha<1$ and therefore $1-\frac{4}{4-\alpha}$ is negative, so we can drop it. Since $q$ and $\alpha$ are both constants, this shows that for any $x_0\in[0,1]$, $\lambda_{n+1}(x_0)$ decays at the rate $O\left(\textstyle \frac{1}{n^\alpha}\right)$ and therefore limits to zero as $n\to \infty$.
	\end{proof}

	We now have the machinery to state the Chebyshev-Markov-Stieltjes inequalities---the key ingredient to our step-function approximation for $F$---which can be found as Theorem 2.5.4 in \cite{Akhiezer1965ClassicalMoment} and Theorem 3.41.1 in \cite{Szego1939}:
	
	\begin{proposition}[Chebyshev-Markov-Stieltjes inequalities] \label{prop:CMSineq}
		Let $F$ be the cumulative distribution function associated with the random variable $X$. If $\lambda_n(x)$ and $\xi_1,\ldots, \xi_n$ are defined as per the preceding discussion, then for each $1\leq j\leq n$,
		\[
		\lambda_n(\xi_1)+\ldots+\lambda_n(\xi_{j-1}) \leq F(\xi_j) \leq \lambda_n(\xi_1)+\ldots+\lambda_n(\xi_{j}),
		\]
		where $\lambda_n(\xi_1)+\ldots+\lambda_n(\xi_{j-1})  = 0$ when $j=1$.
	\end{proposition}
	
	The Chebyshev-Markov-Stieltjes inequalities bound $F$ at each point $\xi_j$, but we can bound $F(x)$ for all $x\in[0,1]$ by forming appropriate step-function bounds. Consider the $x\in[\xi_j,\xi_{j+1}]$. Since $F$ is a cumulative distribution function, we know $F$ is increasing and therefore 
	\begin{equation} \label{eq:boundingFwithUandL}
		\lambda_n(\xi_1)+\ldots+\lambda_n(\xi_{j-1}) \leq F(\xi_j) \leq F(x)\leq F(\xi_{j+1}) \leq \lambda_n(\xi_1)+\ldots+\lambda_n(\xi_{j}).
	\end{equation}
	We can now form a global lower and upper bounds for $F$ by building appropriate step functions. That is, given any $x\in [0,1]$,
	\[
	\sum_{j=1}^{n-1} \lambda_n(\xi_j) \mathbbm{1}_{[\xi_{j+1},1]}(x) \leq F(x) \leq \lambda_n(\xi_1) + \sum_{j=1}^{n-1} \lambda_n(\xi_{j+1}) \mathbbm{1}_{(\xi_{j},1]}(x).
	\]
	In fact, we can show that these bounds converge uniformly to $F(x)$ from below and above, respectively, as $n$ tends to infinity (i.e., as the number of moments considered increases to infinity). Indeed, since $\xi_j\in (0,1)$ for all $j$, Proposition \ref{prop:lambdatozero} guarantees that $\lambda_n(\xi_j)\to 0$ as $n\to \infty$. 
	Thus, with
	\[
	L_n(x) = \sum_{j=1}^{n-1} \lambda_n(\xi_j) \mathbbm{1}_{[\xi_{j+1},1]}(x) \hspace{0.25in} \text{and} \hspace{0.25in} U_n(x) = \lambda_n(\xi_1) + \sum_{j=1}^{n-1} \lambda_n(\xi_{j+1}) \mathbbm{1}_{(\xi_{j},1]}(x)
	\]
	denoting the lower and upper bounds, respectively, we see that 
	\[
	\|U_n(x)-L_n(x)\|_\infty \leq \max_{1\leq j \leq n-1}\{\lambda_n(\xi_j)+\lambda_n(\xi_{j+1})\} \leq 2\max_{1\leq j \leq n}\{\lambda_n(\xi_j)\}\to 0
	\]
	as $n\to \infty$. Therefore, for example, $\|F(x)-L_n(x)\|_\infty \leq \|U_n(x)-L_n(x)\|_\infty \to 0$ as $n\to \infty$ so that the lower bound, $L_n(x)$, converges uniformly to $F(x)$ as $n\to \infty$. Similarly for $U_n(x)$. All together, we've established the following:
	
	\begin{theorem} \label{thm: moment bounds}
		Let $F$ be the (continuous) cumulative distribution function of an arbitrary random variable $X$ defined on $[0,1]$ with moments $\{\mathbb{E}[X^k]\}_{k=0}^{2n-1}$ and let $\{P_i\}_{i=0}^n$ be the associated unique orthonormal polynomials.  If $\xi_1,\ldots,\xi_n$ are the distinct zeros of $P_n$, then
		\[
		\sum_{j=1}^{n-1} \lambda_n(\xi_j) \mathbbm{1}_{[\xi_{j+1},1]}(x) \leq F(x) \leq \lambda_n(\xi_1) + \sum_{j=1}^{n-1} \lambda_n(\xi_{j+1}) \mathbbm{1}_{(\xi_{j},1]}(x),
		\]
		where $\lambda_n(x)= \frac{1}{ \sum_{i=0}^{n-1} |P_i(x)|^2  }$ is the Christoffel function. Moreover, these bounds converge uniformly to $F$ from below and above, respectively, as $n\to \infty$ (i.e., as the number of moments considered increases).
	\end{theorem} 
	
	Proposition \ref{prop:lambdatozero}'s explicit upper bound on the decay rate of $\lambda_n$ allows us to say more. Specifically, with $\alpha= \frac{\ln(p_{\text{max}})}{  \ln\left(\frac{1}{q}\right)}$ being the sharp H\"{o}lder exponent from Theorem \ref{thm:holder1}, we can now provide the exact number of moments needed to approximate $F$ within a given error:
	
	\begin{corollary} \label{corr:momentstoboundF}
		Given $\varepsilon>0$, put $n=\left\lfloor \left(  \frac{8\cdot 16^\alpha (q+1)}{\varepsilon (4-\alpha)}\right)^{1/\alpha} \right\rfloor +2$, where $\alpha= \frac{\ln(p_{\text{max}})}{  \ln\left(\frac{1}{q}\right)}$ is the sharp H\"{o}lder exponent from Theorem \ref{thm:holder1}. If $L_n$ and $U_n$ are computed using the moments $\mathbb{E}[X^k]$ for $0\leq k\leq 2n-1$, then $\|U_n-L_n\|_\infty<\varepsilon$. In particular, so are $\|U_n-F(x)\|_\infty$ and $\|L_n-F(x)\|_\infty$.
	\end{corollary}
	\begin{proof}
		Let $\varepsilon>0$ be given. By Proposition \ref{prop:lambdatozero},
		\begin{align*}
			\|U_n(x)-L_n(x)\|_\infty &\leq \max_{1\leq j \leq n-1}\{\lambda_n(\xi_j)+\lambda_n(\xi_{j+1})\} \\
			&\leq 2\max_{1\leq j \leq n}\{\lambda_n(\xi_j)\} \leq \frac{8\cdot 16^\alpha (q+1)}{(n-1)^\alpha(4-\alpha)},
		\end{align*}
		where $\alpha = \frac{\ln(p_{\text{max}})}{  \ln\left(\frac{1}{q}\right)}$. Thus, the choice
		\[
		n= \left\lfloor \left(  \frac{8\cdot 16^\alpha (q+1)}{\varepsilon (4-\alpha)}\right)^{1/\alpha} \right\rfloor +2,
		\]
		ensures that $\|U_n(x)-F(x)\|_\infty < \varepsilon$. Therefore, as noted in Theorem \ref{thm: moment bounds}, we need to compute the moments $\mathbb{E}[X^k]$ for $0\leq k \leq 2n-1$. The final claim follows from noting that both $\|U_n-F(x)\|_\infty$ and $\|L_n-F(x)\|_\infty$ are less than or equal to $\|U_n-L_n\|_\infty$. 
	\end{proof}

	\bibliographystyle{alpha}
	\bibliography{flip}
	
	\Addresses
	
\end{document}